\documentclass[letterpaper,11pt]{article}

\usepackage[margin=1in]{geometry}
\usepackage[T1]{fontenc}
\usepackage[latin1]{inputenc}
\usepackage[english]{babel}
\usepackage{amsmath}
\usepackage{amssymb}
\usepackage{amsfonts}
\usepackage{array}
\usepackage{float}
\usepackage{graphicx}
\usepackage[font=small,labelfont=bf]{caption}
\usepackage{subcaption}
\usepackage[titletoc,title]{appendix} 
\usepackage{xfrac}
\usepackage{placeins}
\usepackage{bm}
\usepackage{amsthm}
\usepackage[x11names]{xcolor}
\usepackage{diagbox}
\usepackage{arydshln}

\usepackage[hidelinks]{hyperref}
\hypersetup{colorlinks=false}

\numberwithin{equation}{section}

\graphicspath{{figures/}}

\newcommand{\bb}[1]{\mathbb{#1}}
\newcommand{\mc}[1]{\mathcal{#1}}
\newcommand{\mf}[1]{\mathfrak{#1}}

\newcommand{\trial}{\mathcal{U}}

\newcommand{\test}{\mathcal{V}}
\newcommand{\Rv}{\mc{R}_\test}
\newcommand{\Rvr}{\mc{R}_{\test^r}}

\newcommand{\norm}[1]{\left\lVert#1\right\rVert}
\newcommand{\eps}{\varepsilon}
\newcommand{\p}{\partial}
\newcommand{\bs}{\boldsymbol}

\newcommand{\TE}{\text{TE}}

\newcommand{\LP}{\text{LP}}

\def\curl{\nabla \times}

\def\hcurl{\nabla_{\hspace{-1pt}h} \times}
\def\hHcurl{H(\text{curl}, \Omega_h)}
\def\Hcurl{H(\text{curl}, \Omega)}
\def\lb{\langle}
\def\rb{\rangle}
\DeclareMathOperator*{\argmin}{arg\,min}

  \definecolor{ICES}{RGB}{94, 156, 174}
  \definecolor{ORANGE}{RGB}{191, 87, 0}
  \definecolor{RED}{RGB}{190, 30, 49}
  \definecolor{SUN}{RGB}{227, 81, 51}
  \definecolor{GREEN}{RGB}{0, 171, 86}
  \definecolor{BLUE}{RGB}{11, 78, 179}
  \definecolor{BROWN}{RGB}{122, 80, 40}
  \definecolor{GREY}{RGB}{50, 50, 50}
\let\oldparagraph=\paragraph
\renewcommand\paragraph[1]{\oldparagraph{#1.}}

\title{A numerical study of the pollution error and DPG adaptivity for long waveguide simulations}

\author{Stefan Henneking and Leszek Demkowicz\\
\small Oden Institute for Computational Engineering and Sciences,\\[-3pt]
\small The University of Texas at Austin, 204 E 24th St, Austin TX 78712, USA,\\[-3pt]
\small E-mail: \href{mailto:stefan@oden.utexas.edu}{stefan@oden.utexas.edu}, \href{mailto:leszek@oden.utexas.edu}{leszek@oden.utexas.edu}}

\date{\today}

\begin{document}

\maketitle

\begin{abstract}
\noindent
High-frequency wave propagation has many important applications in acoustics, elastodynamics, and electromagnetics. Unfortunately, the finite element discretization for these problems suffers from significant numerical pollution errors that increase with the wavenumber. It is critical to control these errors to obtain a stable and accurate method. We study the effect of pollution for very long waveguide problems in the context of robust discontinuous Petrov--Galerkin (DPG) finite element discretizations. Our numerical experiments show that the pollution primarily has a diffusive effect causing energy loss in the DPG method while phase errors appear less significant. We report results for 3D vectorial time-harmonic Maxwell problems in waveguides with more than 8000 wavelengths. Our results corroborate previous analysis for the Galerkin discretization of the Helmholtz operator by Melenk and Sauter (2011). Additionally, we discuss adaptive refinement strategies for multi-mode fiber waveguides where the propagating transverse modes must be resolved sufficiently. Our study shows the applicability of the DPG error indicator to this class of problems. Finally, we illustrate the importance of load balancing in these simulations for distributed-memory parallel computing.
\end{abstract}

\section{Introduction}
%
%
\paragraph{Motivation}
It is well-known that an accurate numerical solution for wave problems with high frequency is difficult to obtain. High-frequency acoustic or electromagnetic wave propagation problems in the time-harmonic setting lead to partial differential equations (PDEs) with an indefinite Helmholtz or Maxwell operator, respectively. Therefore, many advanced solver techniques for Hermitian positive definite discrete systems are not directly applicable to these problems. Additionally, finite element discretizations are typically dependent on satisfying the Nyquist stability criterion, implying that all propagating wave frequencies must be ``resolved'' to a certain extent in order to have a stable discretization. Even though advanced finite element methods such as the discontinuous Petrov--Galerkin (DPG) method can circumvent the stability problem and deliver a robust discretization for any wave number \cite{demkowicz2012wavenumber,petrides2019phd,petrides2017multigrid}, they do not eliminate the numerical pollution error in multiple dimensions. Pollution can manifest itself in different forms: commonly, we observe a diffusive effect causing wave attenuation and/or a dispersive effect resulting in a phase shift. It is critical to control the pollution error for obtaining accurate results.

\paragraph{Background and literature}
Numerical pollution in wave propagation has been studied extensively for the Bubnov--Galerkin finite element method, as well as discontinuous Galerkin methods, least-squares methods, and various other approaches (see \cite{babuska1997pollution,ihlenburg1997finite,babuska1999dispersion,ainsworth2006dispersive,melenk2011conv,demkowicz2012wavenumber} and references therein); however, to the best of our knowledge, numerical studies have mostly been limited to acoustic wave problems with a moderate number of wavelengths. In this paper, we are reporting numerical results for the DPG method applied to the 3D vectorial time-harmonic Maxwell problem in long waveguides with more than 8000 wavelengths and high order of approximation. We attempt to corroborate theoretical results by Melenk and Sauter \cite{melenk2011conv} and provide guidance as to how to best discretize wave problems with high frequency.

\paragraph{Outline}
First, we introduce the DPG methodology and we discuss how the pollution error enters the DPG formulation as a perturbation parameter in the best approximation, in contrast to the classical Bubnov--Galerkin method where the perturbation parameter is present in the stability constant. This is one of the key points for using the DPG method in high-frequency wave propagation since it implies that the DPG discretization is stable for any wavenumber. We briefly go over the mathematical setting for the time-harmonic Maxwell equations in a linear waveguide in the context of the broken ultraweak variational formulation. For the numerical pollution study, we are reporting the experiments for a simple case with one propagating mode in a rectangular waveguide with the intent of eliminating other effects unrelated to pollution. We discuss our findings in the light of theoretical pollution error estimates from the literature. For waveguides with multiple propagating modes, we discuss different adaptive refinement strategies; in particular, we show numerical results for multi-mode fiber simulations illustrating the need for mesh adaptivity in the numerical solution to such problems. In the case of a distributed-memory parallel simulation, we emphasize the importance of dynamic load balancing for this particular problem.

\section{The DPG methodology}
%
%
The DPG method of Demkowicz and Gopalakrishnan \cite{demkowicz2011part2,demkowicz2012part3,demkowicz2017dpg} has been used to solve applications in viscoelasticity \cite{fuentes2017viscoelasticity}, acoustic and electromagnetic wave propagation \cite{petrides2017multigrid,nagaraj2018raman,henneking2019fiber}, compressible fluid dynamics \cite{chan2014dpg}, and linear elasticity \cite{keith2016elasticity}. Its stability properties make it particularly applicable to high-frequency wave propagation problems, where pre-asymptotic stability is essential for driving efficient $hp$-adaptivity \cite{petrides2019phd}. Through on-the-fly computation of problem-dependent optimal test functions, the DPG method guarantees a stable discretization for any well-posed linear variational problem.

\subsection{The ideal DPG method}
Consider an abstract variational problem of the form,
\begin{equation}
	\left\{
		\begin{array}{l}
			u \in \trial, \\
			b(u,v) = l(v),\ v \in \test,
		\end{array}
	\right.
	\label{eq:cont-abstract-variational}
\end{equation}
where
$\trial$ (trial space) and $\test$ (test space) are Hilbert spaces, and $b(\cdot,\cdot)$ is a continuous bilinear (sesquilinear) form on $\trial \times \test$ (with continuity constant $M$),
\begin{equation}
	|b(u,v)| \leq M \norm{u}_\trial \|v\|_\test,
\end{equation}
that satisfies the continuous inf--sup condition (with inf--sup constant $\gamma$),
\begin{equation}
	\inf_{\|u\|_\trial = 1} \sup_{\|v\|_\test=1} |b(u,v)| =: \gamma > 0 .
	\label{eq:cont-inf-sup}
\end{equation}
And, the continuous linear (antilinear) form $l(\cdot)$ satisfies the compatibility condition,
\begin{equation}
	l(v)=0\ \forall v \in \test_0, \quad \text{where} \quad 
	\test_0 := \left\{ v \in \test:\ b(u,v)=0\ \forall u \in \trial \right\}.
\end{equation}
Let $\trial'$ and $\test'$ denote the space of continuous linear (antilinear) functionals on $\trial$ and $\test$, respectively. By the Babu\v{s}ka-Ne\v{c}as theorem, the variational problem (\ref{eq:cont-abstract-variational}) is well-posed, i.e., there exists a unique solution $u \in \trial$ that depends continuously upon the data,
\begin{equation}
	\norm{u}_\trial \leq \frac{1}{\gamma} \norm{l}_{\test'} .
\end{equation}

Consider finite-dimensional subspaces $\trial_h \subset \trial$ and $\test_h \subset \test$, where $\dim(\trial_h) = \dim(\test_h)$, and the corresponding discrete abstract variational problem,
\begin{equation}
	\left\{
		\begin{array}{l}
			u_h \in \trial_h \subset \trial, \\
			b(u_h,v_h) = l(v_h),\ v_h \in \test_h \subset \test .
		\end{array}
	\right.
	\label{eq:discrete-abstract-variational}
\end{equation}
If the discrete inf--sup condition is satisfied, i.e.,
\begin{equation}
	\inf_{\|u_h\|_\trial = 1} \sup_{\|v_h\|_\test=1} |b(u,v)| =: \gamma_h > 0 ,
	\label{eq:discrete-inf-sup}
\end{equation}
then the discrete problem (\ref{eq:discrete-abstract-variational}) is well-posed; and by Babu\v{s}ka's theorem \cite{babuska1971error},
\begin{equation}
	\| u - u_h \|_\trial \leq 
	\frac{M}{\gamma_h} \inf_{\omega_h \in \trial_h} \| u - \omega_h \|_\trial ,
\end{equation}
where $u$ is the exact solution of (\ref{eq:cont-abstract-variational}), $M/\gamma_h$ is the stability constant, and the best approximation error ${\inf_{\omega_h \in \trial_h} \| u - \omega_h \|_\trial}$ is measured in the trial norm $\norm{\cdot}_\trial$. The continuous inf--sup condition (\ref{eq:cont-inf-sup}) does \emph{not} in general imply the discrete inf--sup condition (\ref{eq:discrete-inf-sup}).

In the DPG method, the issue of discrete stability is solved by finding a unique test space, called the \emph{optimal test space} $\test^{opt}$. Given any trial space $\trial_h \subset \trial$, define its optimal test space by,
\begin{equation}
	\test^{opt} := T(\trial_h) ,
\end{equation}
where the \emph{trial-to-test} operator $T: \trial \rightarrow \test$ is defined by,
\begin{equation}
	(Tu_h, v)_\test = b(u_h,v) \quad \forall u_h \in \trial_h, v \in \test .
	\label{eq:trial-to-test}
\end{equation}
For any $u_h \in \trial_h$, equation (\ref{eq:trial-to-test}) uniquely defines $Tu_h$ by the Riesz representation theorem. Let $\mc{B}: \trial \rightarrow \test'$ denote the linear operator induced by the form $b(\cdot,\cdot)$,
\begin{equation}
	\lb \mc{B}u, v \rb_{\test' \times \test} = b(u,v),\ v \in \test ,
\end{equation}
where $\lb \cdot, \cdot \rb_{\test' \times \test}$ denotes the duality pairing on $\test' \times \test$. Then, $T = \Rv^{-1} \mc{B}$, where $\Rv: \test \rightarrow \test'$ is the Riesz map. In other words, for every trial function $u_h$, the trial-to-test operator defines a unique optimal test function, $v_h = \Rv^{-1} \mc{B} u_h$. The optimal test functions realize the supremum in the inf--sup condition. Indeed,
\begin{equation}
	\sup_{0 \neq v \in \test} \frac{| b(u_h,v) |}{\| v \|_\test}
	=
	\sup_{0 \neq v \in \test} \frac{| (Tu_h, v) |}{\| v \|_\test}
	\leq
	\| Tu_h \|_\test
	=
	\frac{| b(u_h, v_h) |}{\| v_h \|_\test} .
\end{equation}
Therefore, $\gamma_h \geq \gamma$, i.e., discrete stability is guaranteed by construction.

\subsection{Breaking the test space}
In the discussion so far, we have neglected computational aspects of the DPG method.
One question that arises immediately in the context of practicality is the cost of the inversion of the global Riesz map $\Rv$. Let $\Omega$ denote the global domain of interest, with Lipschitz boundary $\Gamma \equiv \p \Omega$, and let $\Omega_h$ denote a suitable finite element triangulation of $\Omega$, with mesh skeleton $\Gamma_h$. By ``breaking'' the test space, i.e., employing a larger discontinuous test space, $\test(\Omega_h) \supset \test(\Omega)$, the inversion of the Riesz map on $\Omega$ is localized and can be done independently element-wise. The element-local computational costs are still significant but they can be parallelized efficiently and fast integration techniques can be implemented to accelerate computations by more than one order of magnitude \cite{mora2019fast,badger2019fast}. By reducing the regularity of the test space, new interface terms arise on the mesh skeleton with interface unknowns $\hat u$. The resulting variational problem is,
\begin{equation}
	\left\{
		\begin{array}{l}
			u \in \trial, \hat u \in \hat \trial, \\
			b(u,v) + \lb \hat u, v \rb_{\Gamma_h} = l(v),\ v \in \test(\Omega_h),
		\end{array}
	\right.
	\label{eq:broken-abstract-variational}
\end{equation}
where $\lb \cdot, \cdot \rb_{\Gamma_h}$ denotes an appropriate duality pairing on the mesh skeleton. The new interface unknowns may be interpreted as Lagrange multipliers \cite{demkowicz2017dpg}.

The stability of the formulation with broken test spaces is inherited from the continuous problem. In particular, the broken formulation (\ref{eq:broken-abstract-variational}) is well-posed with a mesh-independent stability constant of the same order as the stability constant for the continuous problem \cite{carstensen2016breaking}.

\subsection{The practical DPG method}
Until now, the trial-to-test operator $T$ has only been defined in the infinite-dimensional setting (\ref{eq:trial-to-test}). To compute optimal test functions in practice, the inversion of the Riesz map must be approximated on a \emph{truncated} finite-dimensional test space $\test^r \subset \test$, $\dim (\test^r) \gg \dim(\trial_h)$, also called the \emph{enriched test space}.

With the Riesz map defined on the truncated test space, $\Rvr: \test^r \rightarrow (\test^r)'$, the approximate trial-to-test operator $T^r: \trial_h \rightarrow \test^r$ is defined by,
\begin{equation}
	T^r := \Rvr^{-1} \iota^T \mc{B},
\end{equation}
where $\iota: \test^r \rightarrow \test$ is the inclusion map. 

Consequently, the practical DPG method with optimal test functions solves,
\begin{equation}
	\left\{
		\begin{array}{l}
			u_h \in \trial_h \subset \trial, \\
			b(u_h,v_h) = l(v_h),\ v_h \in \test^r = T^r \trial_h ,
		\end{array}
	\right.
	\label{eq:practical-abstract-variational}
\end{equation}
with the additional interface term from (\ref{eq:broken-abstract-variational}) when breaking the test space $V^r$.

Because the optimal test functions are approximated, some stability loss is inevitable. Several papers have addressed the issue of controlling and quantifying the stability loss in the DPG method \cite{gopala2014practical,nagaraj2017fortin}.

\FloatBarrier

\section{Pollution study}
%
%
We first introduce the functional setting for the broken ultraweak DPG formulation of the time-harmonic Maxwell equations in a linear waveguide. Then, we discuss how the frequency $\omega$ enters the DPG error analysis, and how an $hp$-refinement strategy can effectively control the pollution error. We report numerical results for the propagation of the fundamental mode in a rectangular waveguide.

\subsection{Function spaces}
For our time-harmonic Maxwell model problem, we define the following standard energy spaces on a bounded domain $\Omega \subset \bb{R}^3$ with Lipschitz boundary $\Gamma \equiv \p \Omega$:
\begin{equation}
\begin{split}
	L^2(\Omega) &:= \{ y: \Omega \rightarrow \bb{C}: \| y \| < \infty  \}, \\
	\Hcurl &:= \{ \bs q: \Omega \rightarrow \bb{C}^3: \bs q \in (L^2(\Omega))^3, \curl \bs q \in (L^2(\Omega))^3 \},
\end{split}
\end{equation}
where $\| \cdot \|$ is the $L^2(\Omega)$ norm.

Suppose $\Omega$ is partitioned into a set $\Omega_h$ of open disjoint elements $\{K\}_{K \in \Omega_h}$ with Lipschitz element boundaries $\{\p K\}_{K \in \Omega_h}$. The broken energy spaces defined on the finite element mesh $\Omega_h$ are:
\begin{equation}
\begin{split}
	L^2(\Omega_h) &:= \{ y \in L^2(\Omega): y|_K \in L^2(K)\ \forall K \in \Omega_h \} = L^2(\Omega), \\
	\hHcurl &:= \{ \bs q \in (L^2(\Omega))^3: \bs q|_K \in H(\text{curl},K)\ \forall K \in \Omega_h \} \supset \Hcurl .
\end{split}
\end{equation}

Additionally, we must define energy spaces for the trace unknowns that arise from breaking the test space. These spaces are defined on the mesh skeleton $\Gamma_h$. In the model problem, we need the following trace space:
\begin{equation}
	H^{-1/2}(\text{curl}, \Gamma_h) :=
	\big\{
	\hat{\bs q} \in \hskip -4pt \prod_{K \in \Omega_h} \hskip -2pt H^{-1/2}(\text{curl}, \p K):
	\exists \bs q \in H(\text{curl}, \Omega):
	\gamma_t(\bs{q}|_K) = \hat{\bs q} 
	\big\} ,
\end{equation}
where the continuous and surjective tangential trace operator is defined element-wise \cite{carstensen2016breaking}:
\begin{equation}
	\gamma_t: H(\text{curl}, \Omega_h) \rightarrow
	\prod_{K \in \Omega_h} \hskip -2pt H^{-1/2}(\text{curl}, \p K) .
\end{equation}

Finally, we introduce the minimum energy extension norm for the trace unknowns:
\begin{equation}
	\| \hat{\bs q} \|_{H^{-1/2}(\text{curl}, \Gamma_h)} 
	:= \inf_{\substack{\bs q \in \Hcurl \\ \gamma_t(\bs q|_K) = \hat{\bs q}}}
	\| \bs q \|_{\Hcurl}\ .
	\label{eq:trace-norm}
\end{equation}

\subsection{Problem formulation}
We consider the linear, isotropic time-harmonic Maxwell problem:
\begin{equation}
	\left\{ \begin{array}{lcll}
		\curl \bs E &\hskip -5pt=&\hskip -5pt -i \omega \mu \bs H & \text{ in } \Omega, \\
		\curl \bs H &\hskip -5pt=&\hskip -5pt i \omega \eps \bs E & \text{ in } \Omega, \\
		\bs n \times \bs E &\hskip -5pt=&\hskip -5pt \bs n \times \bs E_0 & \text{ on } \Gamma,
	\end{array} \right.
	\label{eq:linear-Maxwell}
\end{equation}
where $\bs E$ and $\bs H$ are the complex vector-valued time-harmonic electric and magnetic field, respectively; $\omega$ is the angular wave frequency; $\bs n$ is the outward unit normal; and $\eps$ and $\mu$ are the scalar-valued electric permittivity and magnetic permeability, respectively. And we assume sufficiently regular boundary data. 


\subsubsection{Ultraweak formulation}
The ultraweak variational formulation is obtained by testing (\ref{eq:linear-Maxwell}) with test functions ($\bs F, \bs G$), integrating over $\Omega$, and relaxing both equations:
\begin{equation}
	\left\{ \begin{array}{lcll}
		\multicolumn{4}{l}{\bs E, \bs H \in (L^2(\Omega))^3,}\\
		(\bs E, \curl \bs F) + (i \omega \mu \bs H, \bs F) &\hskip -5pt=&\hskip -5pt -\lb \bs n \times \bs E_0, \bs F \rb_\Gamma , & \bs F \in \Hcurl, \\
		(\bs H, \curl \bs G) - (i \omega \eps \bs E, \bs G) &\hskip -5pt=&\multicolumn{2}{l}{\hskip -5pt 0,\quad \bs G \in \Hcurl: \bs n \times \bs G = 0 \text{ on } \Gamma .}
	\end{array} \right.
	\label{eq:linear-Maxwell-UW}
\end{equation}

\subsubsection{Broken ultraweak formulation}
By breaking the test space, we introduce new trace unknowns on the mesh skeleton $\Gamma_h$. The broken ultraweak variational formulation is:
\begin{equation}
	\left\{ \begin{array}{lcll}
		\multicolumn{4}{l}{\bs E, \bs H \in (L^2(\Omega))^3,}\\
		\multicolumn{4}{l}{\hat{\bs E} \in \hat{\trial}_1,\ \hat{\bs H} \in \hat{\trial}_2,}\\
		(\bs E, \hcurl \bs F) + \lb \bs n \times \hat{\bs E}, \bs F \rb_{\Gamma_h} + (i \omega \mu \bs H, \bs F) &\hskip -5pt=&\hskip -5pt 0, & \bs F \in \hHcurl, \\
		(\bs H, \hcurl \bs G) + \lb \bs n \times \hat{\bs H}, \bs G \rb_{\Gamma_h} - (i \omega \eps \bs E, \bs G) &\hskip -5pt=&\hskip -5pt 0, & \bs G \in \hHcurl ,
	\end{array} \right.
	\label{eq:linear-Maxwell-brokenUW}
\end{equation}
where $h$ denotes element-wise operations, and,
\begin{align}
	\hat \trial_1 &:= \big\{ \hat{\bs q} \in H^{-1/2}(\text{curl}, \Gamma_h): \bs n \times \hat{\bs q} = \bs n \times \bs E_0 \text{ on } \Gamma \big\} , \\
	\hat \trial_2 &:= H^{-1/2}(\text{curl}, \Gamma_h) .
\end{align}

\subsection{Pollution estimates}
The mathematical setting for different variational formulations of the time-harmonic Maxwell equations in context of DPG is analyzed in much detail in \cite{carstensen2016breaking}. We recap a few points that are relevant for our discussion regarding the ultraweak formulation.

Define the following group variables:
\begin{equation}
	\mf u = (\bs E, \bs H),\quad \mf v = (\bs F, \bs G) .
\end{equation}

The Maxwell operator from (\ref{eq:linear-Maxwell}) can be written as:
\begin{equation}
	A \mf u = \left( \curl \bs E + i \omega \mu \bs H,\ \curl \bs H - i \omega \eps \bs E \right) ;
\end{equation}
then, the bilinear form corresponding to the ultraweak formulation (\ref{eq:linear-Maxwell-UW}) is:
\begin{equation}
	b(\mf u, \mf v) = (\mf u, A^* \mf v) ,
\end{equation}
with the formal Adjoint operator defined by:
\begin{equation}
	A^* \mf v = \left( \curl \bs F + i \omega \eps \bs G,\ \curl \bs G - i \omega \mu \bs F \right) .
\end{equation}

In the ideal ultraweak DPG method with unbroken test spaces, the Adjoint graph norm, $\| \mf v \|_\test = \| A^* \mf v \|$, for the test space yields the \emph{optimal test norm}; with this norm, the method delivers the $L^2$ projection. For the broken formulation (\ref{eq:linear-Maxwell-brokenUW}), we define an additional trace group variable, $\hat{\mf u} = (\hat{\bs E}, \hat{\bs H})$, and obtain the bilinear form:
\begin{equation}
	\hat b(\hat{\mf u}, \mf v) =
	\lb \bs n \times \hat{\bs E}, \bs F \rb_{\Gamma_h} + \lb \bs n \times \hat{\bs H}, \bs G \rb_{\Gamma_h} .
\end{equation}
The broken ultraweak variational problem (\ref{eq:linear-Maxwell-brokenUW}) can then be written as:
\begin{equation}
\left\{ \begin{array}{l}
	\mf u \in \trial, \hat{\mf u} \in \hat \trial , \\
	b(\mf u, \mf v) + \hat b(\hat{\mf u}, \mf v)  = l(\mf v),\quad \mf v \in \test ,
\end{array} \right.
\end{equation}
where the load $l(\mf v) \equiv \bs 0$ in this case.

The optimal test norm is not localizable, but we can augment it with an additional term to obtain a quasi-optimal test norm: $\| \mf v \|^2_\test = \| A^* \mf v \|^2 + \alpha \| v \|^2$, with scaling parameter $\alpha \in \mathcal{O}(1)$. The quasi-optimal test norm is robustly equivalent with the optimal test norm, i.e. independent of the frequency $\omega$, and the robust stability constant is maintained in the broken formulation \cite{carstensen2016breaking}. This implies that the approximation error is bounded by the best approximation error (BAE) uniformly in $\omega$:

\begin{equation}
	\| \mf u - \mf u_h \|^2 + \| \hat{\mf u} - \hat{\mf u}_h \|^2_{\hat q} \leq
	C \left[ \inf_{\mf w_h} \| \mf u - \mf w_h \|^2 + \inf_{\hat{\mf w}_h} \| \hat{\mf u} - \hat{\mf w}_h \|^2_{\hat q} \right] ,
	\label{eq:error-bound}
\end{equation}

where constant $C$ is independent of the mesh and frequency $\omega$, and $\| \cdot \|_{\hat q}$ refers to an operator-dependent minimum energy extension norm, related to (\ref{eq:trace-norm}). The estimate also implies that the $L^2$ best approximation is pollution free because it is independent of $\omega$. In one dimension, the BAE for the traces is zero, thus the method is in fact pollution free \cite{demkowicz2012wavenumber,petrides2017multigrid}. In multiple dimensions, however, the BAE for the traces, measured in $\| \cdot \|_{\hat q}$, does depend on the frequency $\omega$ and the method exhibits numerical pollution.

The estimate for the standard Galerkin method on the other hand has a stability constant that is not $\omega$-independent, thus the Galerkin discretization is not robustly stable. The DPG method hides the perturbation parameter $\omega$ in the best approximation and by doing so yields a \emph{stable discretization for any wavenumber}. This can practically be exploited by starting computation on a coarse mesh where the pollution error is high, and driving $hp$-adaptivity with the DPG error indicator. This approach yields superior meshes for resolving localized waves \cite{petrides2017multigrid,petrides2019phd}.

A wavenumber explicit analysis for the Helmholtz equation is presented for the DPG method in \cite{demkowicz2012wavenumber} and for the Galerkin method in \cite{melenk2011conv}. For the Galerkin discretization, Melenk and Sauter show that quasi-optimality is obtained under the conditions that $\omega h/p$ is sufficiently small and $p$ is at least $\mathcal{O}(\log \omega)$, where $h$ is the mesh size and $p$ the polynomial order of approximation \cite{melenk2011conv}. Based on these estimates, the best approach to dealing with the pollution error may be an $hp$-strategy that preferably increases the polynomial order $p$ while keeping $\omega h$ constant for increasing frequency.

In the next section, we study the pollution error with numerical experiments for many wavelengths and discuss the observations with regard to the suggested $hp$-strategy and its applicability to the DPG method for the time-harmonic Maxwell problem.

\subsection{Numerical results}
\begin{figure}[htb]
	\centering
	\begin{subfigure}[b]{0.4\textwidth}
		\includegraphics[width=\textwidth]{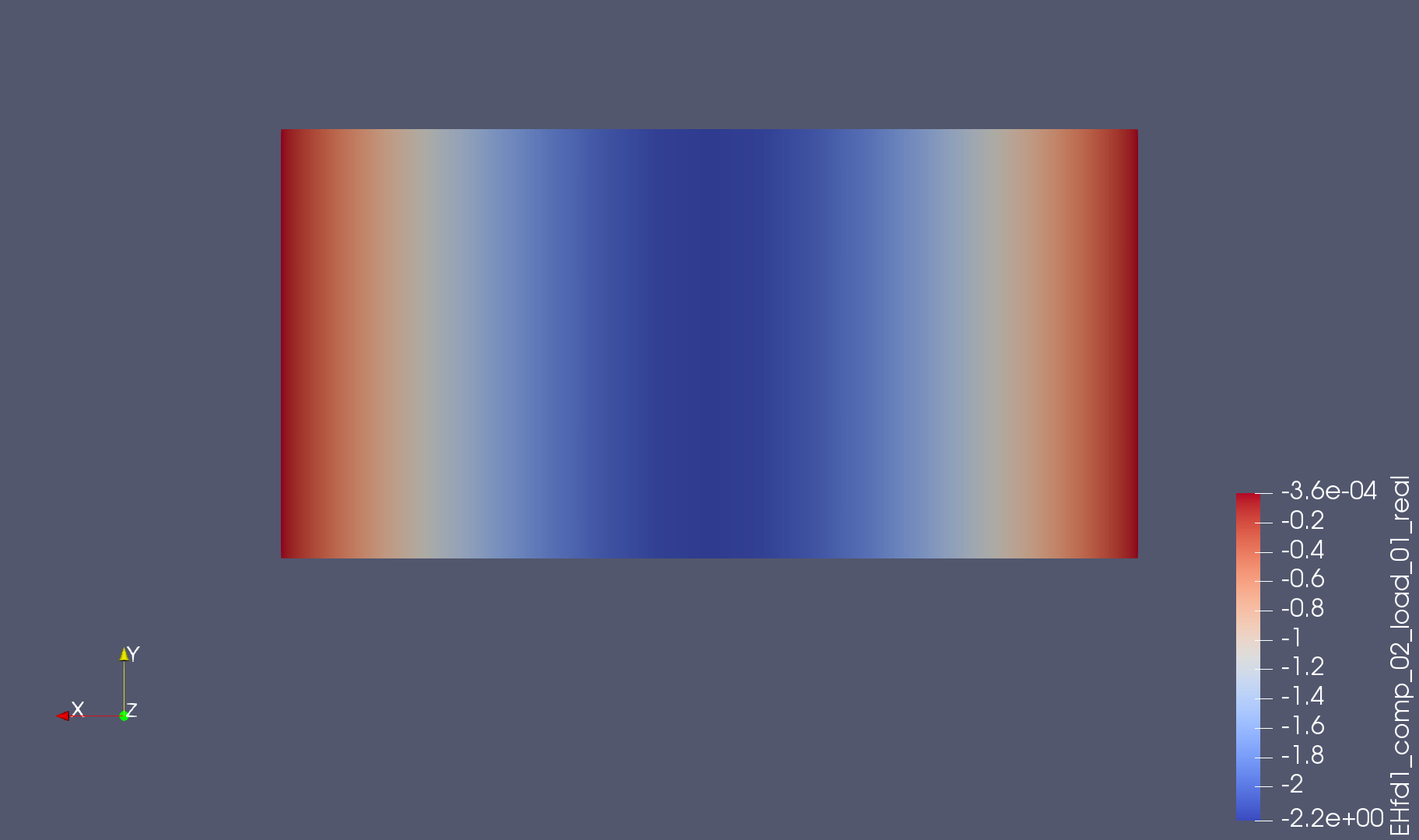}
		\caption{Electric field component $E_y$}
		\label{fig:TE10-E_y}
	\end{subfigure}
	\begin{subfigure}[b]{0.4\textwidth}
		\includegraphics[width=\textwidth]{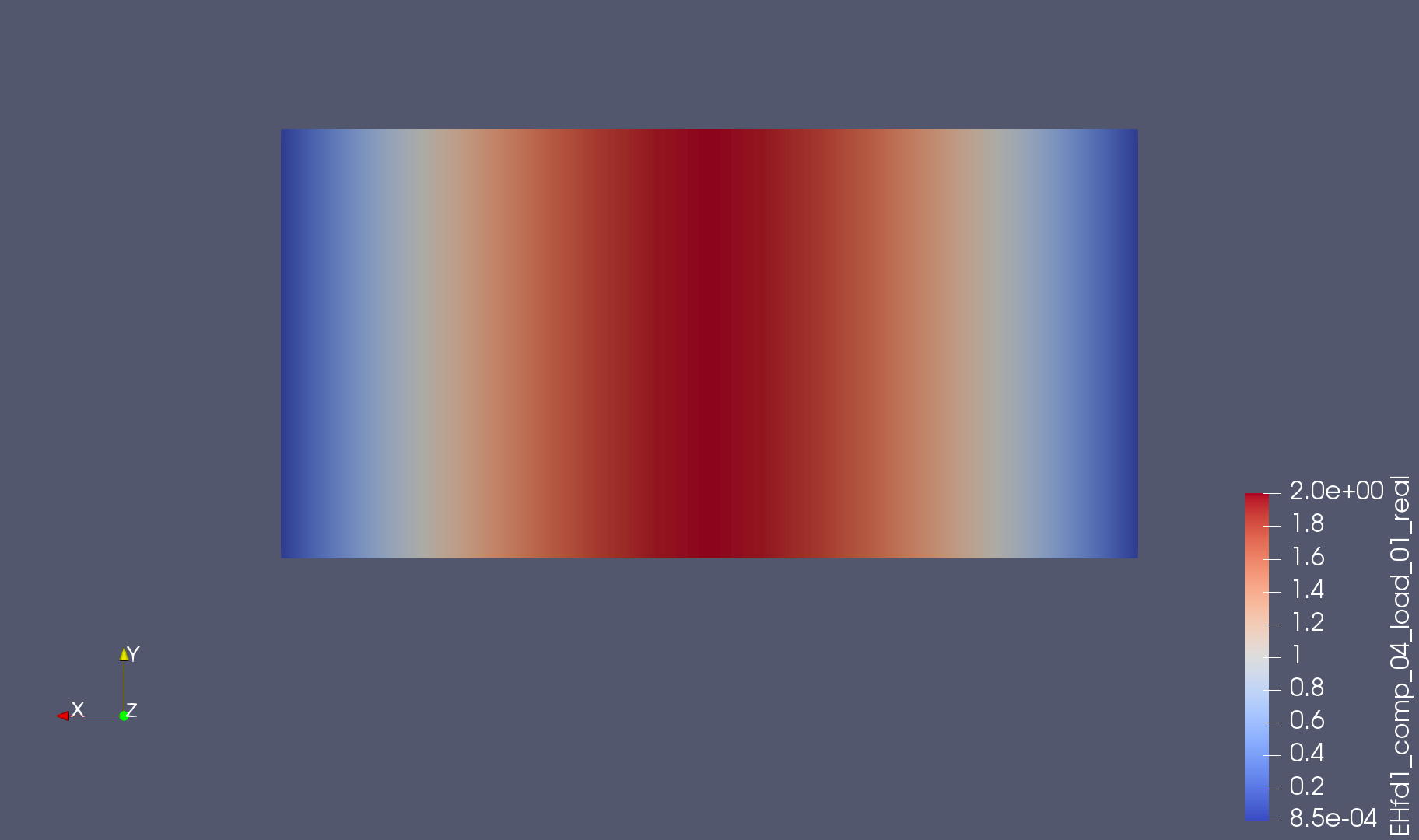}
		\caption{Magnetic field component $H_x$}
		\label{fig:TE10-H_x}
	\end{subfigure}
	\caption{$\text{TE}_{10}$ transverse fields in a rectangular waveguide in a plane normal to the $z$-axis. The simple transverse profile of the fundamental mode justifies a geometry discretization with few elements in the cross-section.}
	\label{fig:TE10-transverse}
\end{figure}

The propagation of an electromagnetic field in a waveguide is governed by the Maxwell equations. We assume that the time-harmonic setting is justified and that the waveguide medium is nonmagnetic, dielectric, and no free charges are present. While nonlinear effects and anisotropic, inhomogeneous material properties play important roles in research on fiber optics, for the purpose of this pollution study we will assume the waveguide medium is linear, isotropic and homogeneous. If we prescribe idealized PEC boundary conditions, then under these assumptions the Maxwell equations reduce to (\ref{eq:linear-Maxwell}). Note that in fact we do prescribe non-zero Dirichlet boundary conditions (BCs) at the input to excite the waveguide and impedance BCs at the output but PEC BCs everywhere else. For this simplified setting, the propagating field in a waveguide can be described as a superposition of guided modes. These modes are eigenfunction solutions to the reduced scalar Helmholtz problem which can be solved analytically for simple domains (e.g, see \cite{griffiths,jackson}). Consider a rectangular waveguide with the domain $\Omega$, in Cartesian coordinates:
\[
	\Omega = (0.0,1.0) \times (0.0,0.5) \times (0,L) ,
\]
where $L$ is the length of the waveguide. The fundamental mode in this waveguide is the transverse electric $\TE_{10}$ mode, depicted in Fig.\ \ref{fig:TE10-transverse}. The fundamental mode is not very oscillatory in the transverse direction. At the waveguide end, we employ an absorbing impedance boundary condition that matches the wave impedance for the fundamental mode. In the rectangular waveguide experiments, the cross-section is modeled with two hexahedral elements, which is justified by the simple transverse mode profile (cf. Fig.\ \ref{fig:TE10-transverse}).

\begin{figure}[htb]
	\centering
	\includegraphics[width=.8 \textwidth]{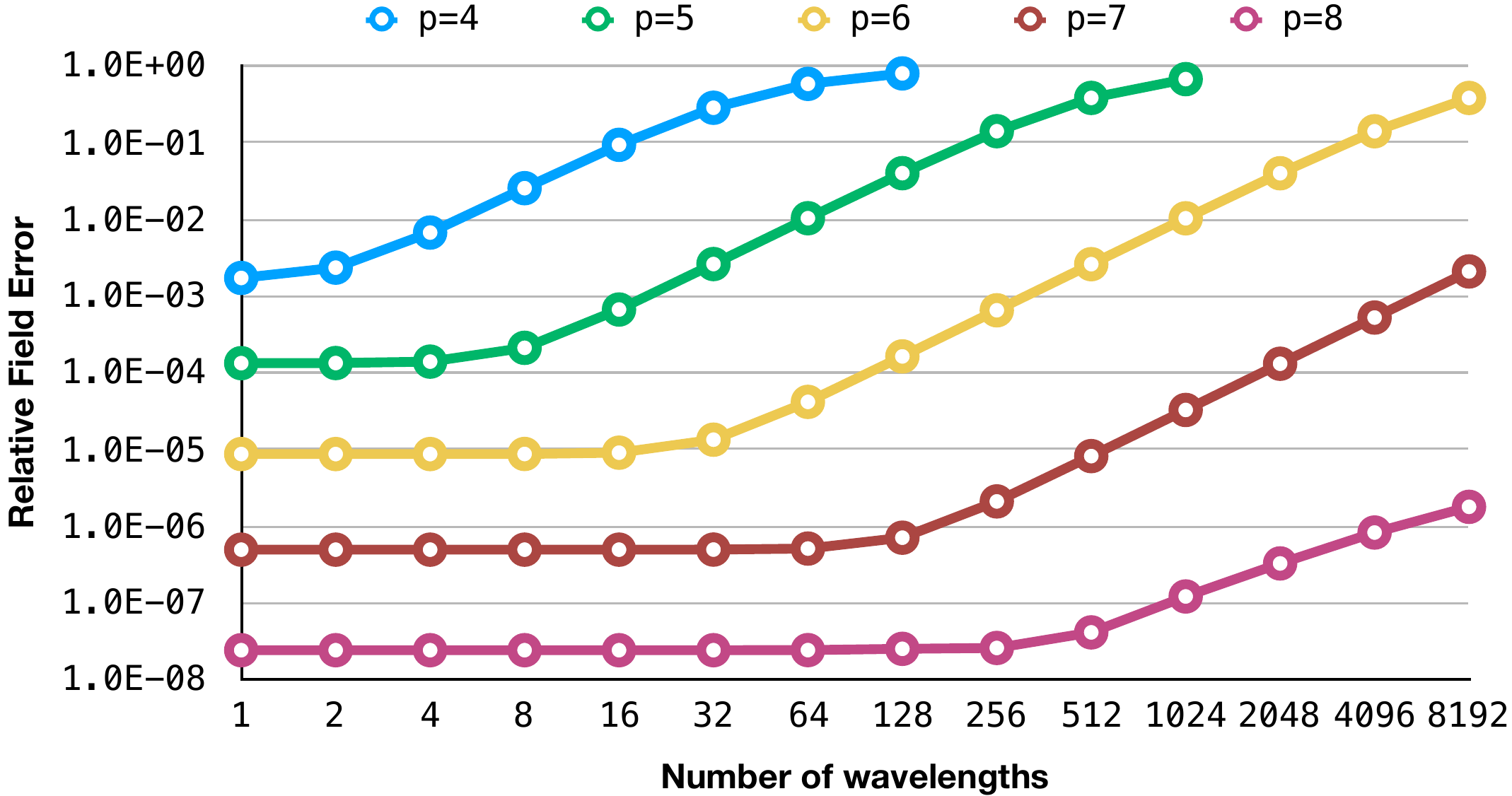}
	\caption{Relative field error for uniform order of approximation. As the number of wavelengths is increased, we keep the ratio of elements per wavelength constant. To counter the pollution error, the polynomial order of approximation must be increased logarithmically with the number of wavelengths.}
	\label{fig:rect-uniform-p-error}
\end{figure}

In our first experiment, we analyze the relative field error, measured in the $L^2$ norm, for the propagating fundamental mode in waveguides of different length $L$. The smallest waveguide has a length equivalent to one wavelength of the fundamental mode, and the longest one has 8192 wavelengths. As we increase the length $L$, we keep the number of elements per wavelength (i.e., degrees of freedom (DOFs) per wavelength) constant. In particular, we choose a discretization with four elements per wavelength. Fig.~\ref{fig:rect-uniform-p-error} shows the relative field error for these waveguides for uniform order of approximation, ranging from $p=4$ to $p=8$. In all numerical experiments, we are using the enrichment order $\Delta p = 1$ for the test space to approximate optimal test functions.

We make several observations: for a fixed number of wavelengths, higher polynomial order yields significantly smaller (more than one order of magnitude) errors, as expected; for every order of approximation, the field error starts to increase if the waveguide is long enough despite keeping the DOFs per wavelength constant; and for higher $p$, this pollution effect is ``kicking in'' at a later point, i.e., more wavelengths can be computed with higher order before the pollution error is measurable. Furthermore, to maintain some desired accuracy, one needs to increase the polynomial order in nearly regular intervals. For example, to achieve $1\%$ accuracy for 4 wavelengths, it is sufficient to use $p=4$; at 64 wavelengths, $p=5$ is needed; with $p=6$, computing up to 1024 wavelengths is feasible with this error margin; and $p=7$ would most likely be sufficient for 16384 wavelengths. At a closer look, these intervals resemble a logarithmic dependency on the polynomial order $(4 * 2^4 = 64,\ 64*2^4 = 1024,\ 1024*2^4=16384)$. In other words, these results corroborate theoretical estimates by Melenk and Sauter predicting that control of the pollution error would require increasing $p$ logarithmically with the wavenumber. In the pollution regime, we observe that the relative error grows as $\mathcal{O}((\#\lambda)^2)$, as expected, for all $p$ except $p=8$ where the rate appears to slow; currently, we do not have an explanation for the latter observation.

\begin{figure}[htb]
	\centering
	\includegraphics[width=.8 \textwidth]{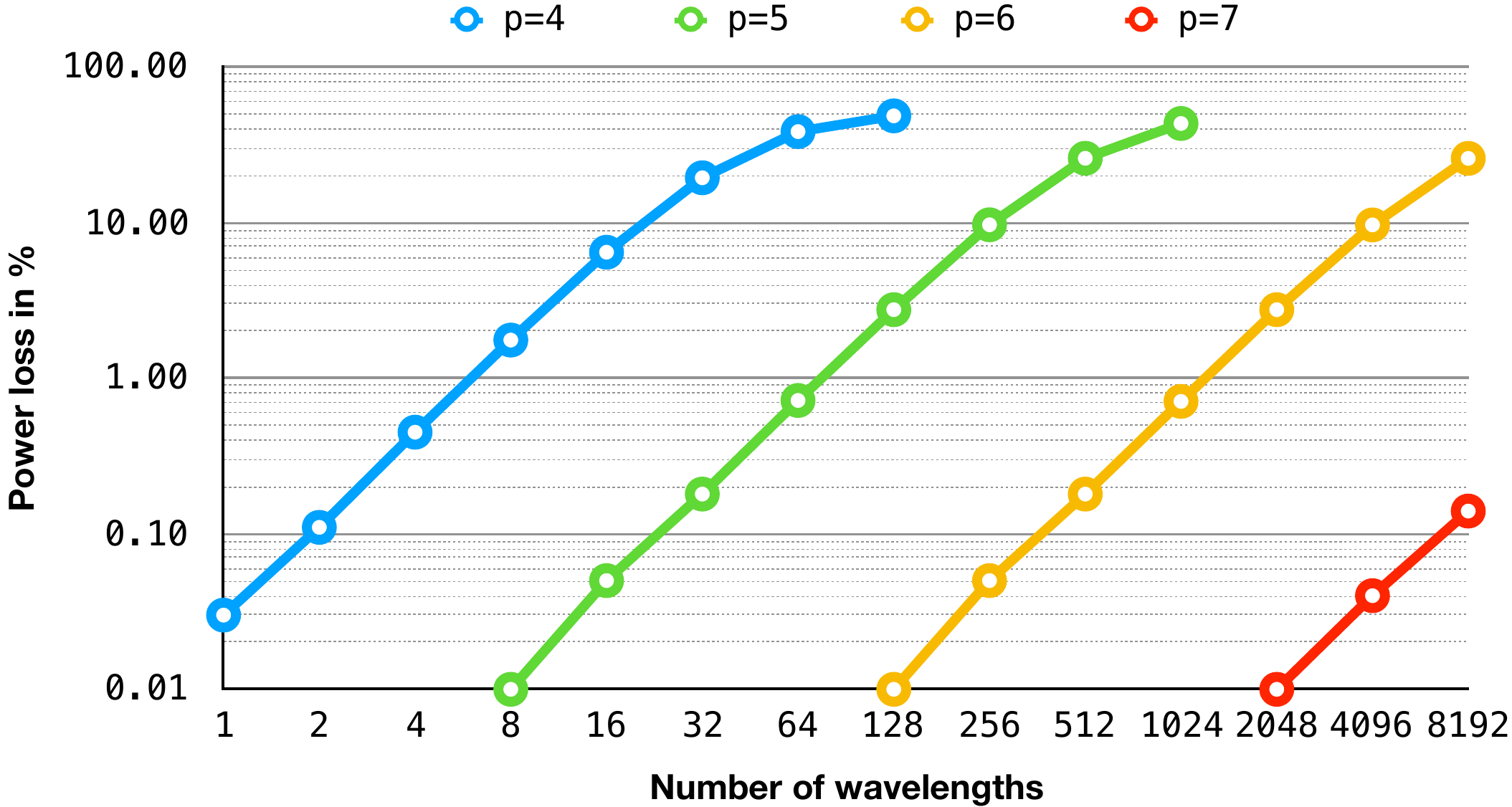}
	\caption{Power loss for uniform order of approximation. In the DPG discretization, we observe that the pollution error has primarily a diffusive effect. This implies that the power flux, measured perpendicular to the waveguide cross-section, decreases along the waveguide. The depicted quantity is the loss of power between the input and output for waveguides of different lengths.}
	\label{fig:rect-uniform-p-power}
\end{figure}

In our experiments, the pollution was primarily a diffusive error causing wave attenuation. This is in agreement with previous observations for the DPG method \cite{demkowicz2012wavenumber}. A practical way of measuring this diffusivity in waveguide applications is to compute power flux through the cross-section of the waveguide at different points in $z$. In a linear, dielectric waveguide with PEC boundary conditions, the fundamental mode should be carried without loss of power. Fig.~\ref{fig:rect-uniform-p-power} shows the measured power loss between the waveguide input $(z=0)$ and output $(z=L)$ for different polynomial orders. Note that $p=8$ has less than $0.005\%$ loss of power in all tested waveguides. The pollution error is clearly visible in terms of power loss. We also observe the same logarithmic dependency for increasing polynomial order, illustrated by the near-equidistant parallel character of the lines.

\begin{figure}[htb]
	\centering
	\begin{subfigure}[b]{0.4\textwidth}
		\includegraphics[width=\textwidth]{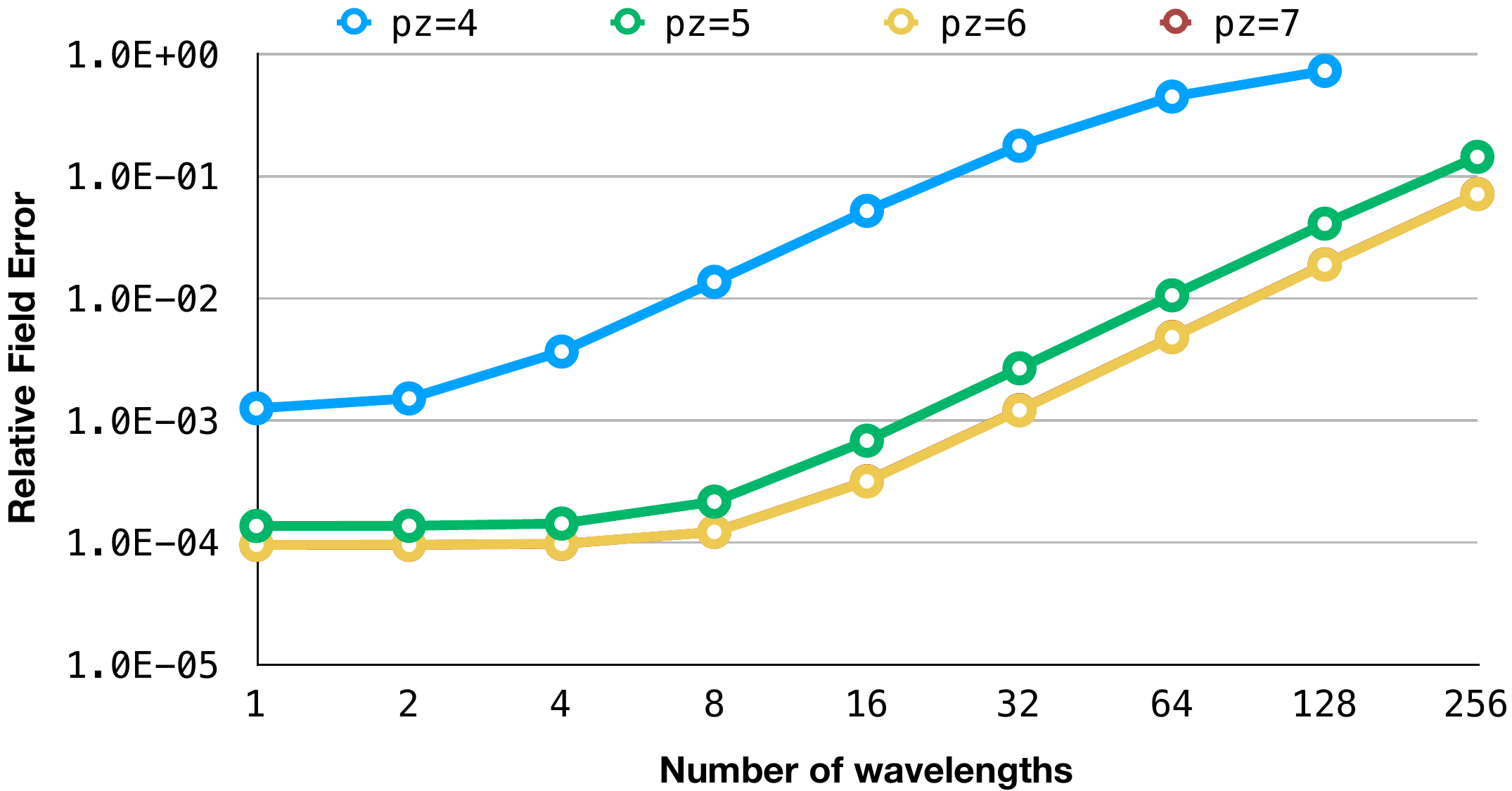}
		\caption{Anisotropic approximation order $p$}
		\label{fig:rect-aniso-p-error}
	\end{subfigure}
	\begin{subfigure}[b]{0.4\textwidth}
		\includegraphics[width=\textwidth]{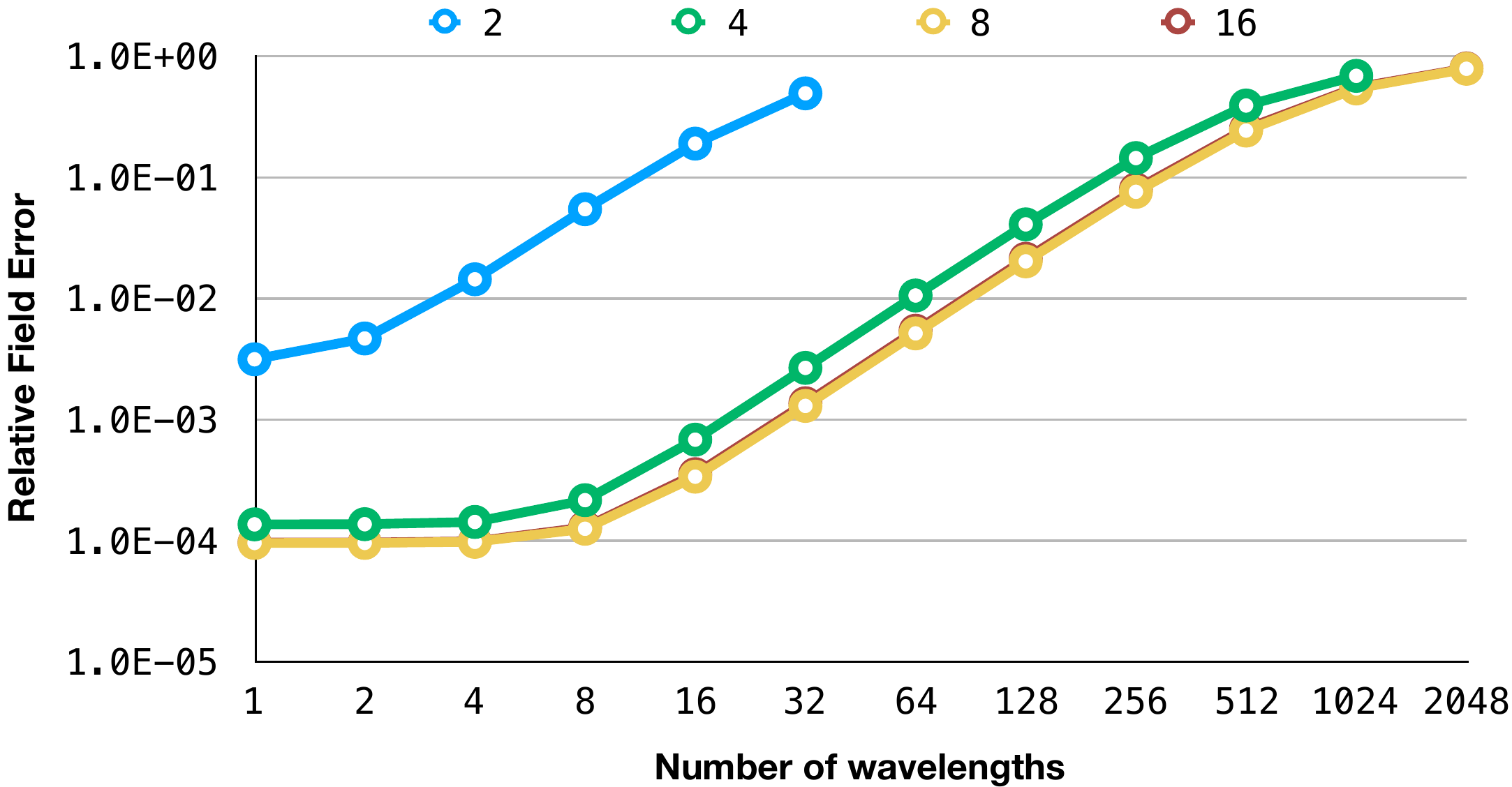}
		\caption{Anisotropic element size $h$}
		\label{fig:rect-aniso-h-error}
	\end{subfigure}
	\caption{Relative field error for anisotropic refinements. Increasing the number of DOFs anisotropically via $h$- or $p$-refinements (in the direction of propagation) does not suffice to control the pollution error in our experiments.}
	\label{fig:rect-aniso-error}
\end{figure}

Moving on to additional experiments, we keep our focus on the same rectangular waveguide but with different potential approaches of dealing with the pollution error. It may be reasonable to assume that since the wave is propagating in one direction (along $z$), it will be sufficient to increase the order of approximation anisotropically or to increase the number of elements through anisotropic $h$-refinements in $z$. Exploring both of these options (cf.\ Fig.\ \ref{fig:rect-aniso-error}), we find that neither one of these approaches yields satisfactory results. First, in Fig.\ \ref{fig:rect-aniso-p-error}, we use fifth-order polynomials in the radial discretization ($px=py=5$) of the waveguide and increase the anisotropic order from $pz=4$ to $pz=7$. While the error decreases initially, it begins stagnating at $pz=6$ (note that the $pz=7$ error coincides almost exactly with $pz=6$). The same observation is made for uniform order $p=5$ with varying number of elements per wavelength (ranging from 2 elements to 16 elements). Our findings indicate that the pollution error comes from an interplay between the mode resolution (radial discretization) and the wave resolution in the direction of propagation. In other words, \emph{increasing the number of DOFs anisotropically does not suffice asymptotically to control the pollution error.}

\FloatBarrier

\begin{figure}[htb]
	\centering
	\begin{subfigure}[b]{0.4\textwidth}
		\includegraphics[width=\textwidth]{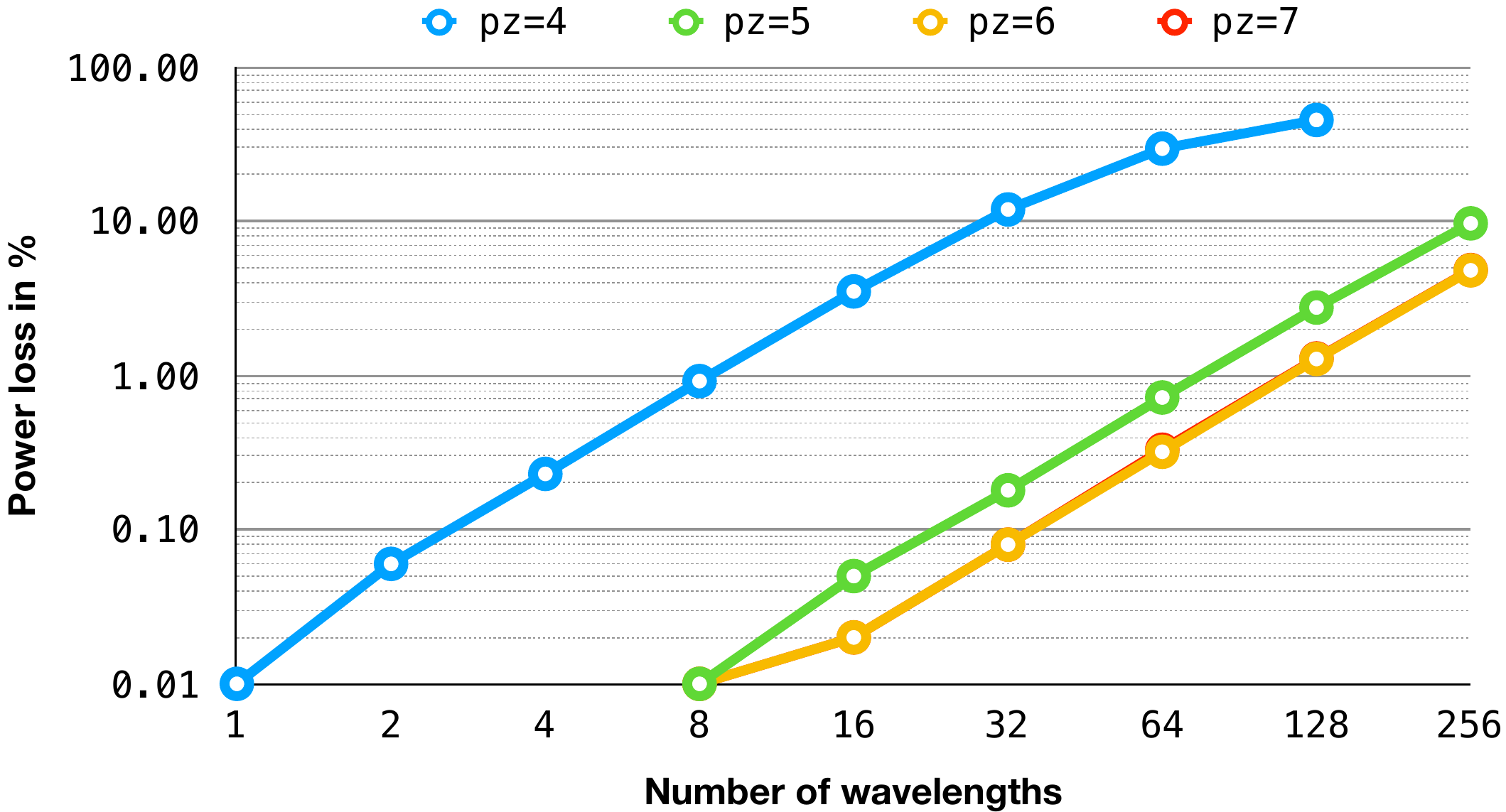}
		\caption{Anisotropic approximation order $p$}
		\label{fig:rect-aniso-p-power}
	\end{subfigure}
	\begin{subfigure}[b]{0.4\textwidth}
		\includegraphics[width=\textwidth]{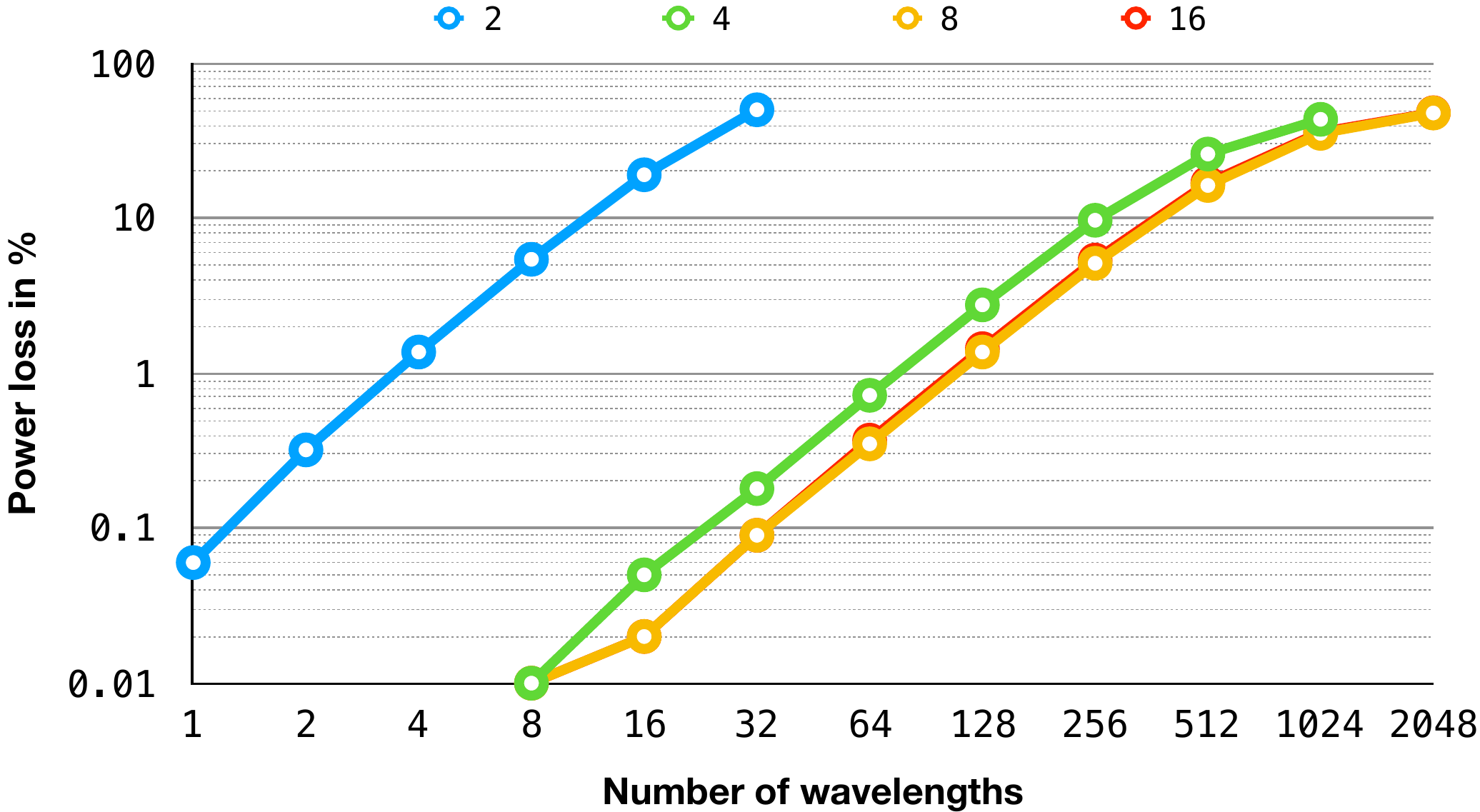}
		\caption{Anisotropic element size $h$}
		\label{fig:rect-aniso-h-power}
	\end{subfigure}
	\caption{Power loss for anisotropic refinements. The diffusive pollution effect is not countered through anisotropic $h$- or $p$-refinements.}
	\label{fig:rect-aniso-power}
\end{figure}

Finally, we measure the loss of power for both anisotropic refinement cases, plotted in Fig.\ \ref{fig:rect-aniso-power}. Expectably, we observe the same stagnation in the diffusive pollution, consistent with the errors measured in the previous plot.

We have conducted these experiments on different waveguides (rectangular waveguides, circular waveguides, and step-index fibers) with various propagating modes, and the observations are all consistent with the observations presented up to here; we therefore omit showing additional numerical results for those cases.

\vskip -10pt

\FloatBarrier

\section{Adaptivity study}
%
%

In our adaptivity study, we focus on a different aspect of resolving the propagating wave. We have shown that the interplay between resolving the wave along the direction of propagation and resolving the transverse mode profile is important in controlling the pollution error. For that reason, the finite element mesh should be sensitive to different mode profiles and adapt to resolve them appropriately. This is especially important in waveguide applications where significant transfer of power occurs between different guided modes. Our tool for adapting the mesh ``on-the-fly'' is the DPG residual that serves as an error estimator in the energy norm. We briefly recap the important points in the derivation of this indicator and proceed with numerical experiments in multi-mode step-index fibers. Lastly, we will look at the load imbalance that results from adapting the mesh to different propagating modes.

\subsection{The DPG error indicator}

The DPG method can be reformulated as a minimum residual method \cite{demkowicz2017dpg},
\begin{equation}
		u_h 
		= \argmin_{w_h \in \trial_h} \| l - \mc{B} w_h \|_{\test'}^2
		= \argmin_{w_h \in \trial_h} \| \Rv^{-1} (l - \mc{B} w_h) \|_{\test}^2\ ,
	\label{eq:abstract-min-residual}
\end{equation}
where the residual is minimized in the dual test norm $\| \cdot \|_{\test'}$. Taking the G\^{a}teaux derivative, we obtain a minimum residual formulation,
\begin{equation}
	\left\{ \begin{array}{l}
		u_h \in \trial_h, \\
		( \Rv^{-1} (l-\mc{B}u_h), \Rv^{-1} \mc{B} w_h )_\test = 0, \quad w_h \in \trial_h .
	\end{array} \right.
\end{equation}
Furthermore, we define the energy norm $\norm{\cdot}_E$ on the trial space $\trial$ by,
\begin{equation}
	\| u \|_E := \| \mc{B}u \|_{\test'} = \| \Rv^{-1} \mc{B}u \|_{\test} .
	\label{eq:energy-norm}
\end{equation}
We define $\psi$ as the Riesz representation of the residual,
\begin{equation}
	\psi := \mc{R}_\test^{-1} (l - \mc{B} u_h) .
	\label{eq:error-representation}
\end{equation}
Notice that when $u_h$ minimizes the residual (\ref{eq:abstract-min-residual}), then,
\begin{equation}
	(\psi, \Rv^{-1} \mc{B} w_h )_\test = 0, \quad w_h \in \trial_h .
\end{equation}
We arrive at a mixed Galerkin formulation,
\begin{equation}
\arraycolsep=2pt
	\left\{ \begin{array}{lll}
		u_h \in \trial_h,\ \psi \in \test, \\
		(\psi, v)_\test + b(u_h, v) &= l(v), &\quad v \in \test, \\
		b(w_h, \psi) &= 0, &\quad w_h \in \trial_h .
	\end{array} \right.
\end{equation}
The error measured in the energy norm can be computed explicitly,
\begin{equation}
		\| u - u_h \|_E = \| \mc{B}(u-u_h) \|_{\test'}
	= \| l - \mc{B}u_h \|_{\test'}
	= \| \Rv^{-1} (l - \mc{B}u_h) \|_\test
	= \| \psi \|_\test ,
	\label{eq:DPG-residual}
\end{equation}
hence $\norm{\psi}$ offers a built-in a-posteriori error indicator. Finally, note that the choice of the test norm $\norm{\cdot}_\test$ is critical, as it dictates the norm in which the method converges. The natural choice for the test norm in the ultraweak formulation is the adjoint graph norm \cite{demkowicz2017dpg}.

\subsection{Multi-mode step-index fibers}
\begin{figure}[htb]
	\centering
	\includegraphics[width=0.8 \textwidth]{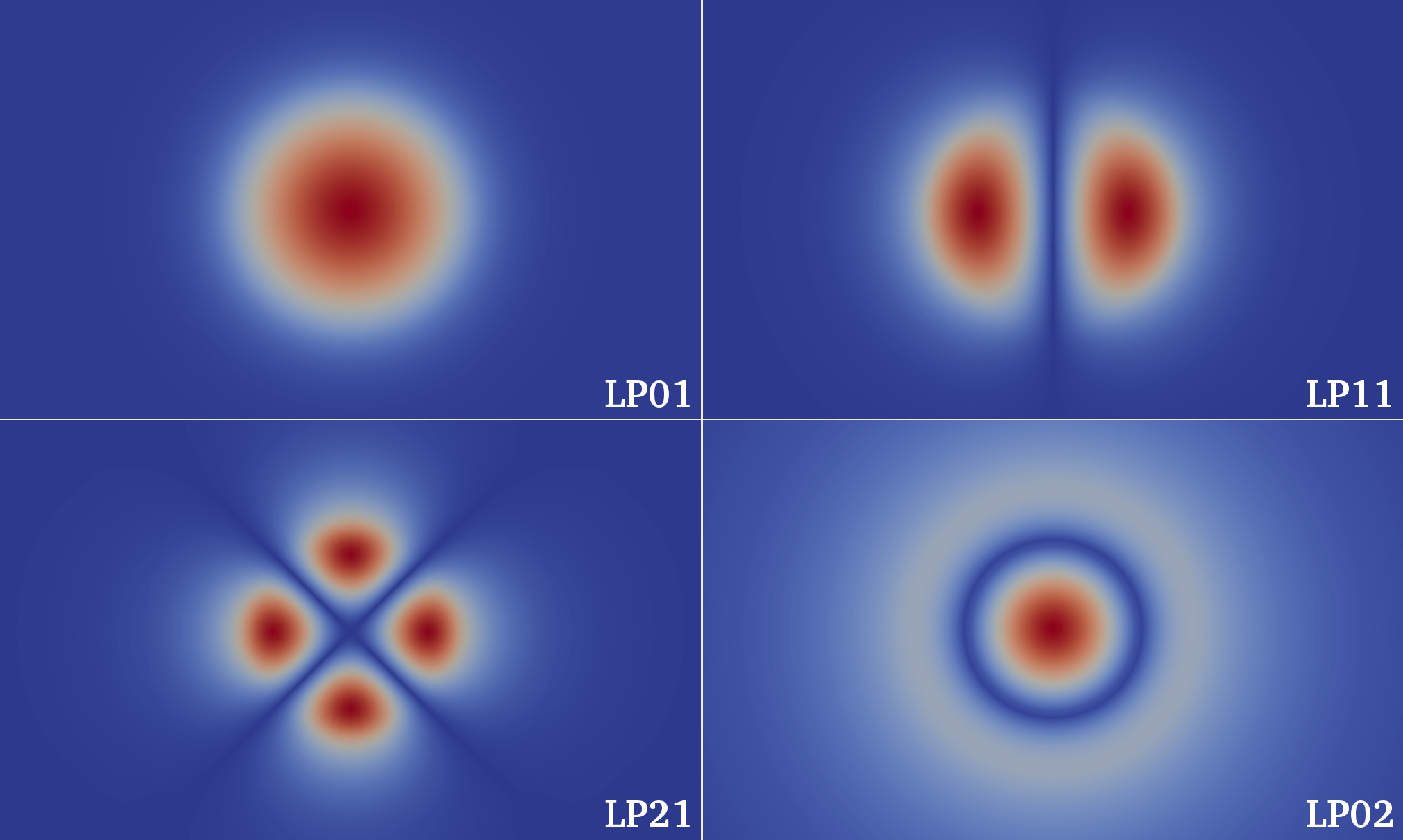}
	\caption{Guided modes in LMA fiber (magnitude of the electric field). Higher-order modes are more oscillatory in the transverse direction and carry more energy outside the fiber core. Therefore, they require additional refinements of the fiber cross-section discretization.}
	\label{fig:modes}
\end{figure}

We consider a dielectric step-index waveguide. More precisely, consider a weakly-guiding large mode area (LMA) step-index fiber made of silica glass. See Tab.\ \ref{tab:fiber-parameters} for a description of the model parameters for the fiber. For weakly-guiding fibers, $(n_{core} - n_{clad}) / n_{core} \ll 1$, and the guided modes are linearly polarized (LP) modes. The $V$-number is a ``normalized frequency'' that determines how many guided modes are supported by the particular fiber. For example, if $V < 2.405$, then the fiber is single-mode. LMA fibers have a relatively large core radius and support multiple modes; the fiber we are using has $V \approx 4.43$, and it supports four guided modes: $\{ \LP_{01},\LP_{11},\LP_{21},\LP_{02} \}$. The fiber axis is assumed to be aligned with the $z$-axis, and the length of the fiber is $L$.  Fig.\ \ref{fig:modes} illustrates the guided modes for this particular fiber, showing the magnitude of the electric field in the center of the fiber cross-section.

For multi-mode propagation, we are using an absorbing perfectly matched layer (PML) at the fiber end. We refer to \cite{astaneh2018pml} for details on the derivation and implementation of a PML for the DPG method.

\begin{table}[htb]
	\centering
	\caption{Step-index fiber: parameters}
	\begin{tabular}{l|l|l}
	Description & Symbol & Value \\
	\hline
	Wavelength & $\lambda$ & $1064$ nm \\
	Core radius & $r_{core}$ & $12.7\ \mu\text{m}$ \\
	Cladding radius & $r_{clad}$ & $127\ \mu\text{m}$ \\
	Core refractive index & $n_{core}$  & $1.4512$ \\
	Cladding refractive index & $n_{clad}$  & $1.45$ \\
	Numerical aperture & NA & $0.059$  \\
	V-number & $V$ & $4.43$
	\end{tabular}
	\label{tab:fiber-parameters}
\end{table}

For this particular fiber, the power confinement (amount of energy confined to the core region) of each mode is,
\begin{equation}
	\{ \Gamma_{01}, \Gamma_{11}, \Gamma_{21}, \Gamma_{02} \} \approx 
	\{ 96.11, 88.77, 74.79, 59.58 \} \%.
\end{equation}

Clearly, the optimal discretization of the fiber cross-section is different for each of these modes. That is, to capture the oscillations of the higher-order modes near the core-cladding interface, a finer discretization is needed than for the fundamental mode $\LP_{01}$. In particular, higher-order modes demand more refinements (or degrees of freedom) outside of the fiber core, when compared to the fundamental mode that is mostly confined to the core region. Therefore, one geometry cannot be optimal for capturing any propagating mode.


Suppose we are interested in simulating the transverse mode instability (TMI) phenomenon \cite{eidam2011experimental} in active gain fiber amplifiers. The TMI is characterized by the chaotic transfer of energy between the fundamental mode and the higher-order modes. One challenge in computing a numerical solution to the resulting nonlinear Maxwell problem is to capture modes accurately when they occur. With mode instabilities, it is not known a-priori which modes will be propagating in which parts of the fiber. Refining the initial geometry globally to better resolve higher-order modes increases the computational cost dramatically and may render the computation infeasible for large problem instances. Adaptivity, on the other hand, can be used to refine the mesh where it is needed for capturing these modes locally, and the overall computational cost will be kept significantly lower.

\FloatBarrier
\subsection{Numerical experiments}
In the following experiments, we are aiming to establish the efficacy of adaptivity based on the DPG residual for resolving higher-order modes. In the broken DPG setting, the residual is computed through element-wise contributions, i.e.,
\begin{equation}
	\| \psi \|_\test^2 = \sum_{j=1}^{n} \| \left. \psi \right|_{K_j}\|_{\test(K_j)}^2\ ,
\end{equation}
where $K_j, j=1,\ldots,n$, denotes the $j$-th element. After each solve, elements are marked for refinement if they satisfy a certain criterion. We use a strategy for marking elements that is based on D\"orfler's marking \cite{dorfler1996marking}:
\begin{enumerate}
	\item{Sort the element residuals $\| \left. \psi \right|_{K_j}\|_{\test(K_j)}^2$ in descending order, i.e., 
\begin{equation}
	\| \left. \psi \right|_{K_1}\|_{\test(K_1)} \geq 
	\| \left. \psi \right|_{K_2}\|_{\test(K_2)} \ldots  \geq 
	\| \left. \psi \right|_{K_n}\|_{\test(K_n)}\ .
\end{equation}
	}
	\item{Mark elements $K_j, j=1,\ldots,J$, where $J \leq n$ is the smallest integer for which the following is true:
\begin{equation}
	\sum_{j=1}^{J} \| \left. \psi \right|_{K_j}\|_{\test(K_j)}^2 \geq \kappa \| \psi \|_\test^2\ ,
	\label{eq:adaptive-marking}
\end{equation}
	where $\kappa \in (0,1)$.}
\end{enumerate}

At this point, with some choice of $\kappa$, elements have been marked for refinement. However, it is not clear how to optimally refine each marked element when the $hp$ mesh supports anisotropic adaptive refinements in both element size $h$ and polynomial order $p$. The choice will ultimately be problem-dependent.

In our experiments, we choose an initial mesh with uniform polynomial order $p=5$, two elements per wavelength in $z$-direction (direction of propagation), and a radial (transverse) hybrid discretization using curvilinear hexahedral and prismatic elements.

\begin{figure}[htb]
	\centering
	\includegraphics[width=0.8 \textwidth]{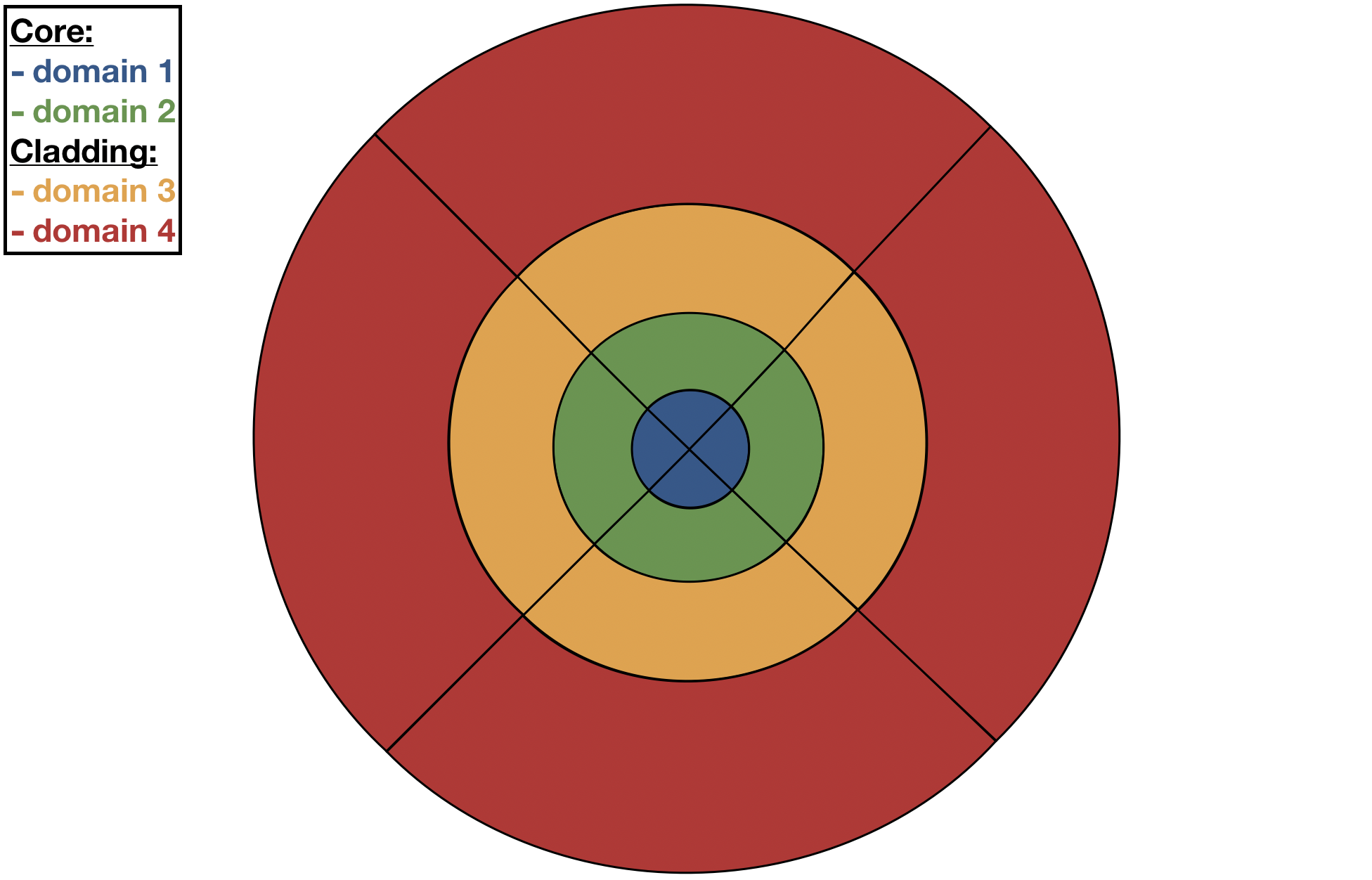}
	\caption{Initial geometry discretization in the fiber cross-section (not drawn to scale). This discretization works well for the fundamental mode, but it may not be sufficient to capture higher-order modes.}
	\label{fig:fiber-domains}
\end{figure}

Fig.\ \ref{fig:fiber-domains} illustrates the initial geometry discretization in the fiber cross-section (not drawn to scale): four prisms are used to model the center of the fiber core, and they are surrounded by three layers of four hexahedra each. We refer to these different layers as ``domains'' and enumerate them from 1 to 4 moving radially outward from the center of the fiber to the cladding boundary. The choice of the initial discretization was informed by the fiber geometry, the fact that all guided modes decay exponentially in the cladding region, and by conducting numerical tests primarily with the fundamental mode. For a relatively short fiber of 16 wavelengths, this initial geometry captures the fundamental mode very well with regard to several physical quantities of interest (e.g., conservation of power, mode confinement). Higher-order modes are not captured quite as well, and for fibers with many wavelengths (i.e., several hundred or a few thousand wavelengths) we observe more significant pollution errors in the propagation of these modes. For example, the errors can be observed in small oscillations of the mode powers along the fiber, diffusive pollution effects, or an unsteady power confinement ratio.

We use DPG to perform multiple adaptive mesh refinements, each based on the respective previous solution and residual, to test the residual error indicator for capturing different modes. As a test case, we will look at the 16 wavelengths fiber. The goal is to observe the sensitivity of the adaptive refinements toward specific propagating modes. We choose the parameter $\kappa=0.5$ in (\ref{eq:adaptive-marking}) and proceed with four adaptive refinement steps after the initial solution. Note that mesh regularity requirements may perform some additional refinements to ``close the mesh'', i.e., to obtain a mesh that is 1-irregular.

\begin{figure}[htb]
	\centering
	\begin{subfigure}[b]{0.75\textwidth}
		\includegraphics[width=\textwidth]{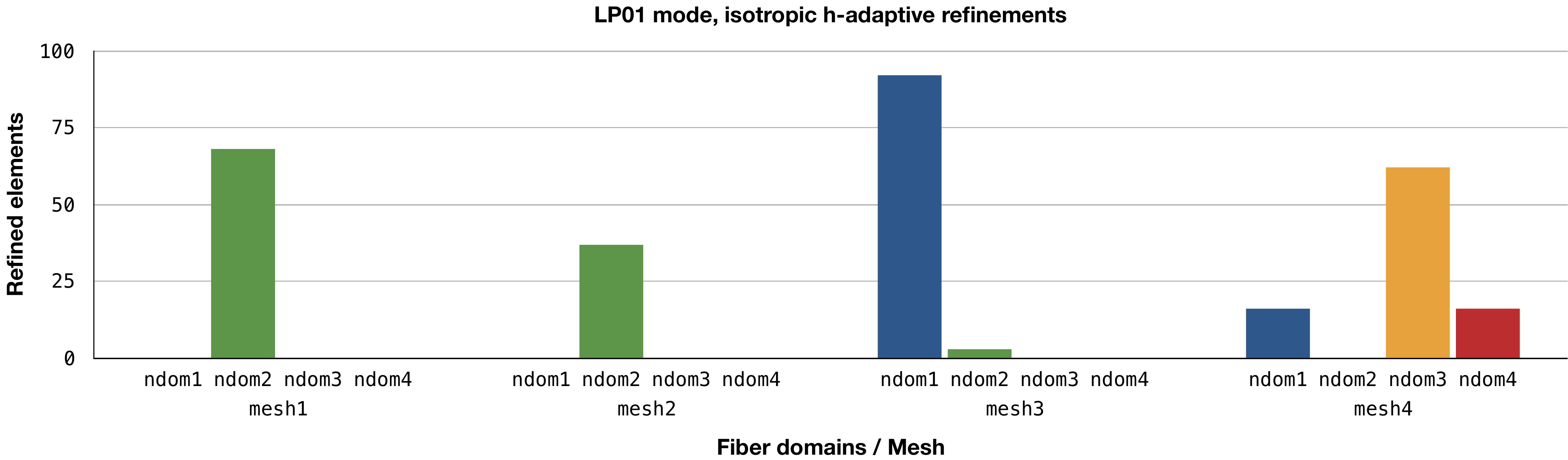}
	\end{subfigure}
	\begin{subfigure}[b]{0.75\textwidth}
		\includegraphics[width=\textwidth]{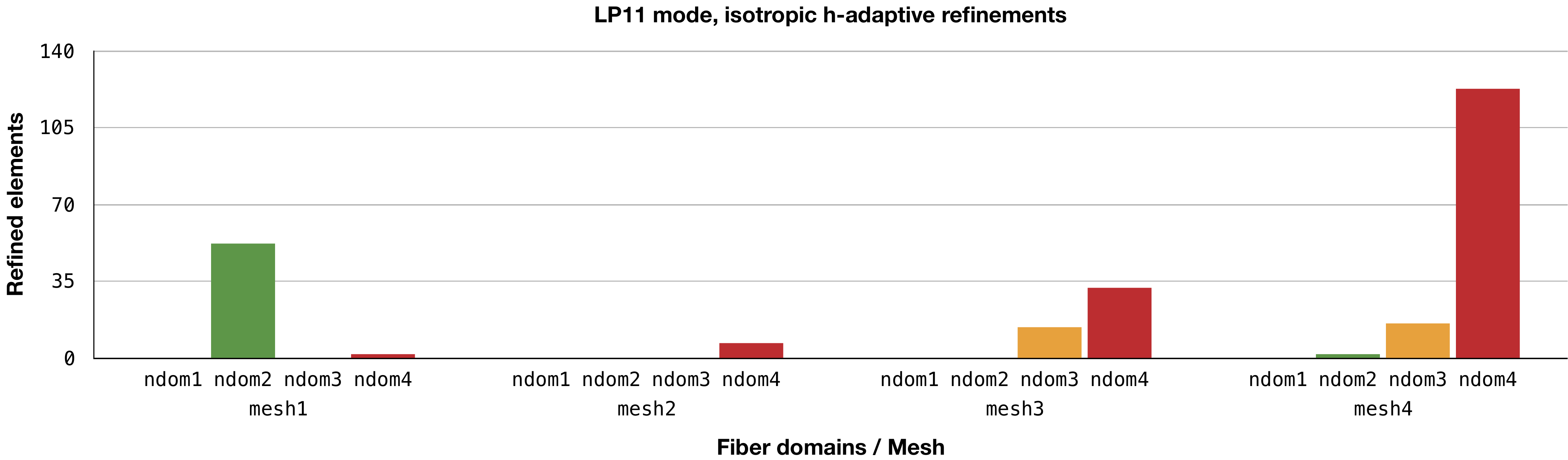}
	\end{subfigure}
	\begin{subfigure}[b]{0.75\textwidth}
		\includegraphics[width=\textwidth]{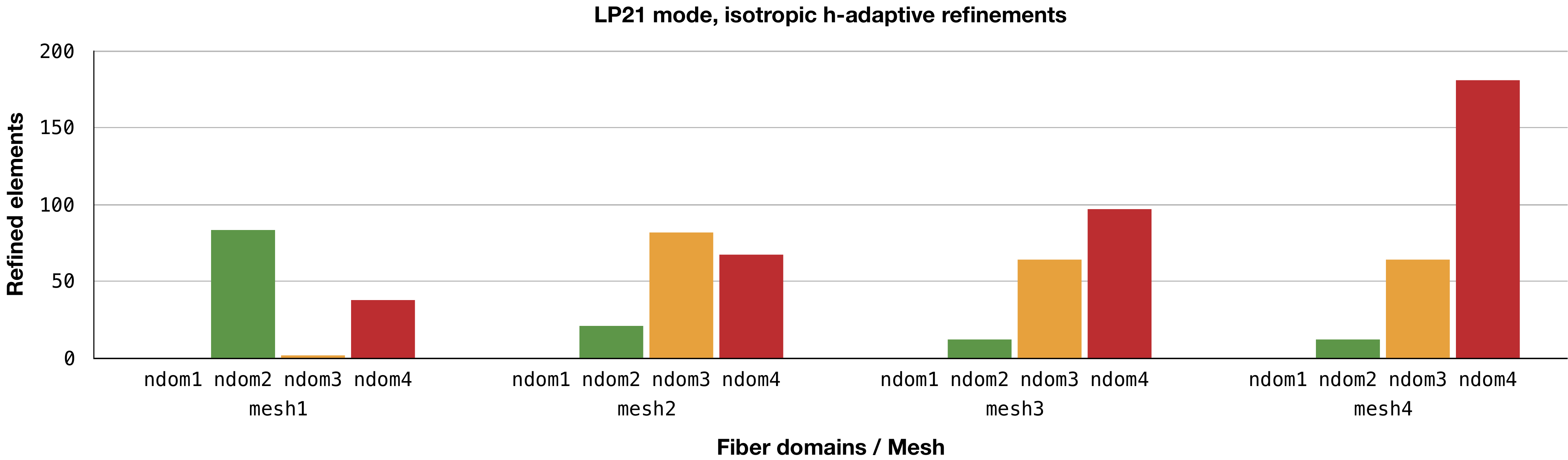}
	\end{subfigure}
	\begin{subfigure}[b]{0.75\textwidth}
		\includegraphics[width=\textwidth]{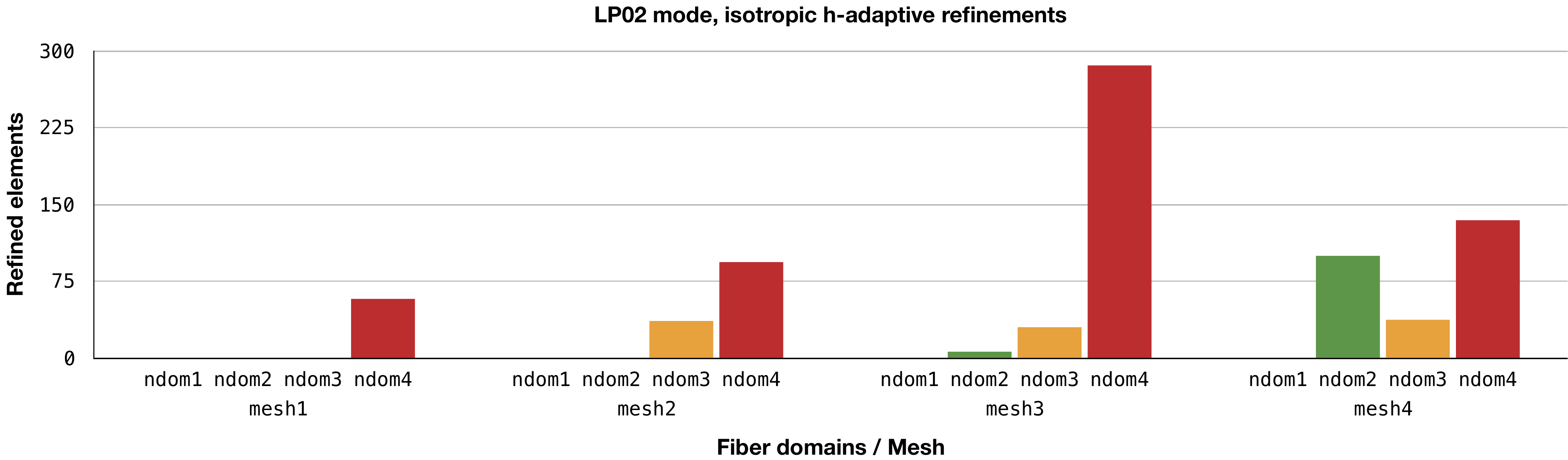}
	\end{subfigure}
	\caption{Isotropic $h$-adaptive refinements. Depending on the propagating mode, the DPG error indicator marks elements for refinement in different fiber domains.}
	\label{fig:iso-adaptive}
\end{figure}

First, we apply isotropic $h$-refinements for marked elements. Fig.\ \ref{fig:iso-adaptive} shows the domains of refinement in the fiber. Each plot illustrates how the mesh is successively refined for one particular guided mode propagating in the fiber. For the fundamental mode, the error indicator is marking elements for refinement primarily in the fiber core, where most of the energy is located. The first three refinement steps exclusively refine in the outer and inner core region. None of the adaptive refinements for higher-order modes refine inside the inner core region. It is notable how sensitive the error indicator is to these different modes. In the case of the $\LP_{02}$ mode, the code primarily refines in the outer cladding region; this is likely to be the case because the initial cladding geometry discretization is too coarse to handle the exponential decay of the remaining energy in the transverse field.

\begin{figure}[htb]
	\centering
	\begin{subfigure}[b]{0.75\textwidth}
		\includegraphics[width=\textwidth]{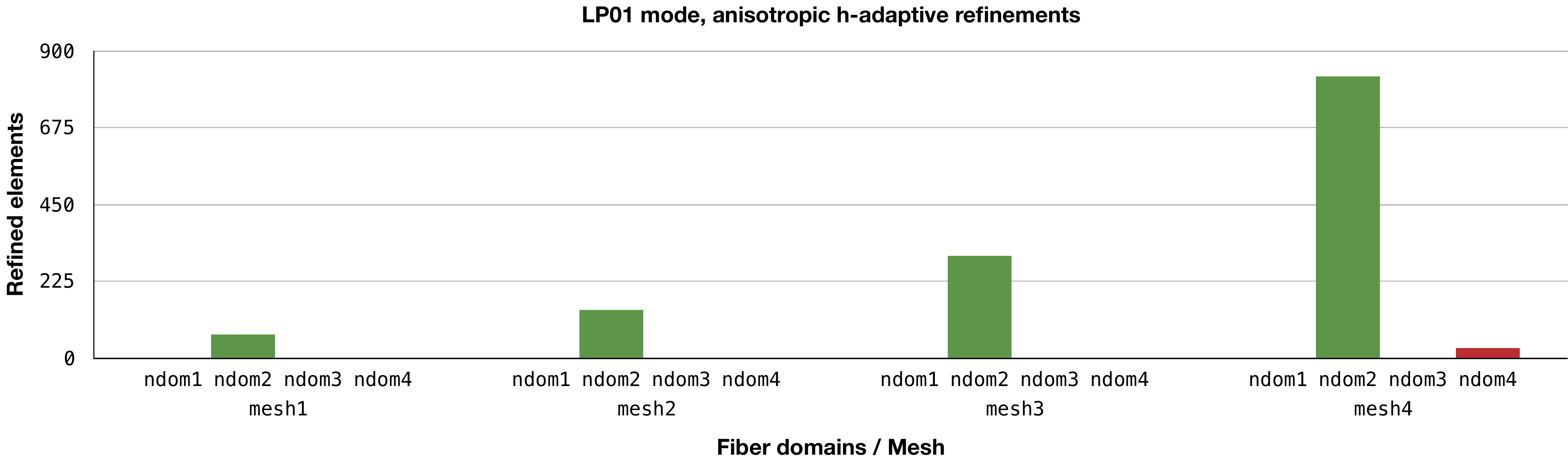}
	\end{subfigure}
	\begin{subfigure}[b]{0.75\textwidth}
		\includegraphics[width=\textwidth]{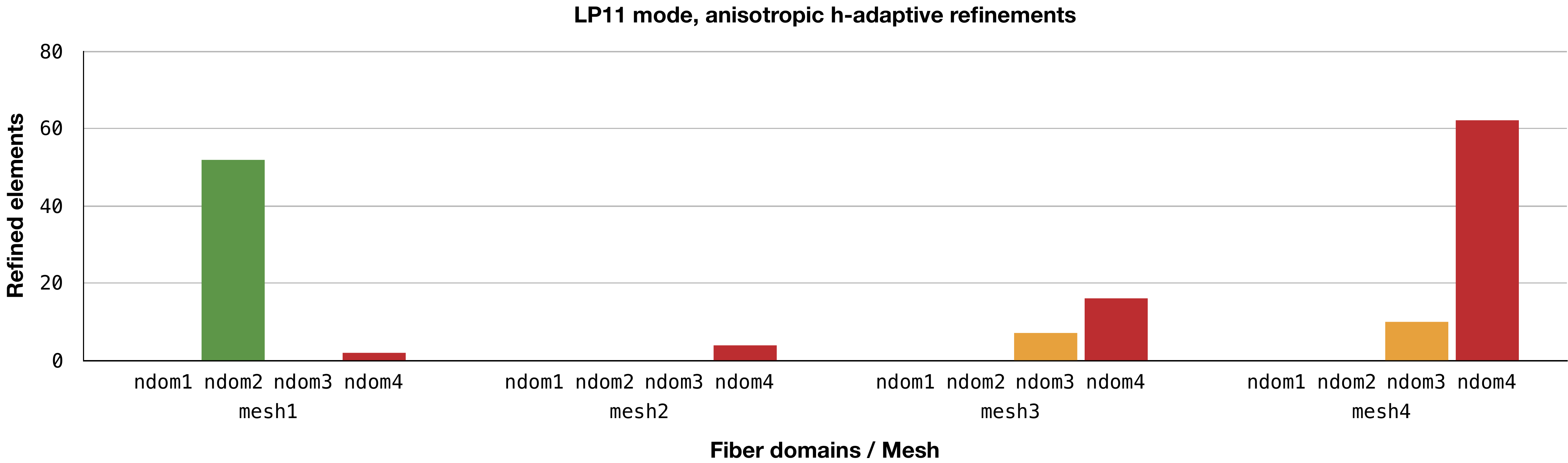}
	\end{subfigure}
	\begin{subfigure}[b]{0.75\textwidth}
		\includegraphics[width=\textwidth]{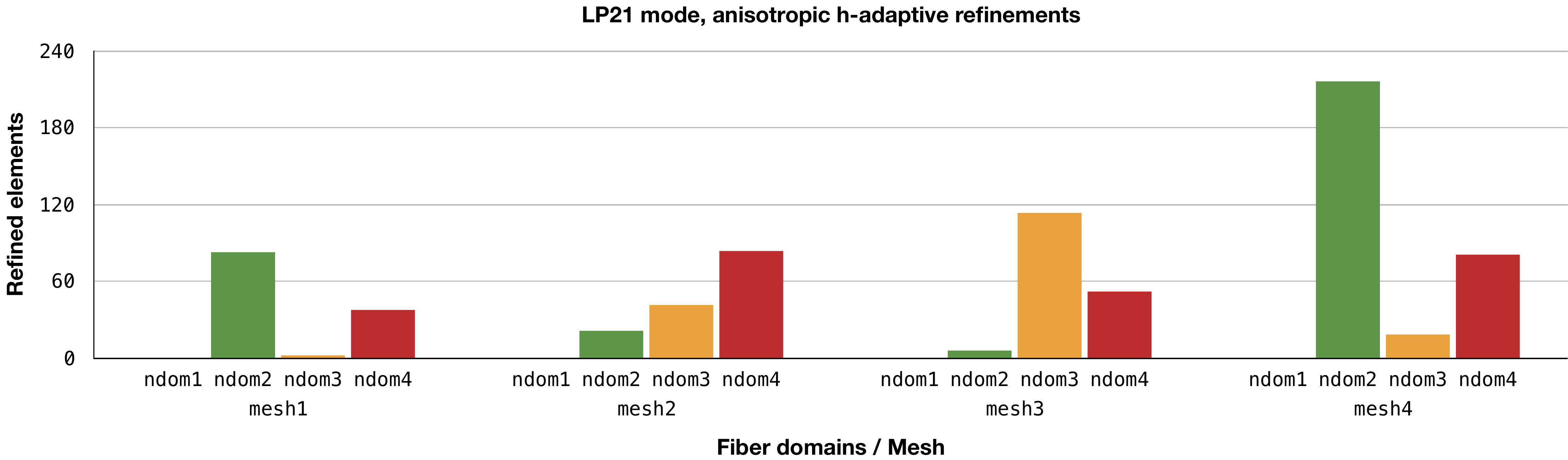}
	\end{subfigure}
	\begin{subfigure}[b]{0.75\textwidth}
		\includegraphics[width=\textwidth]{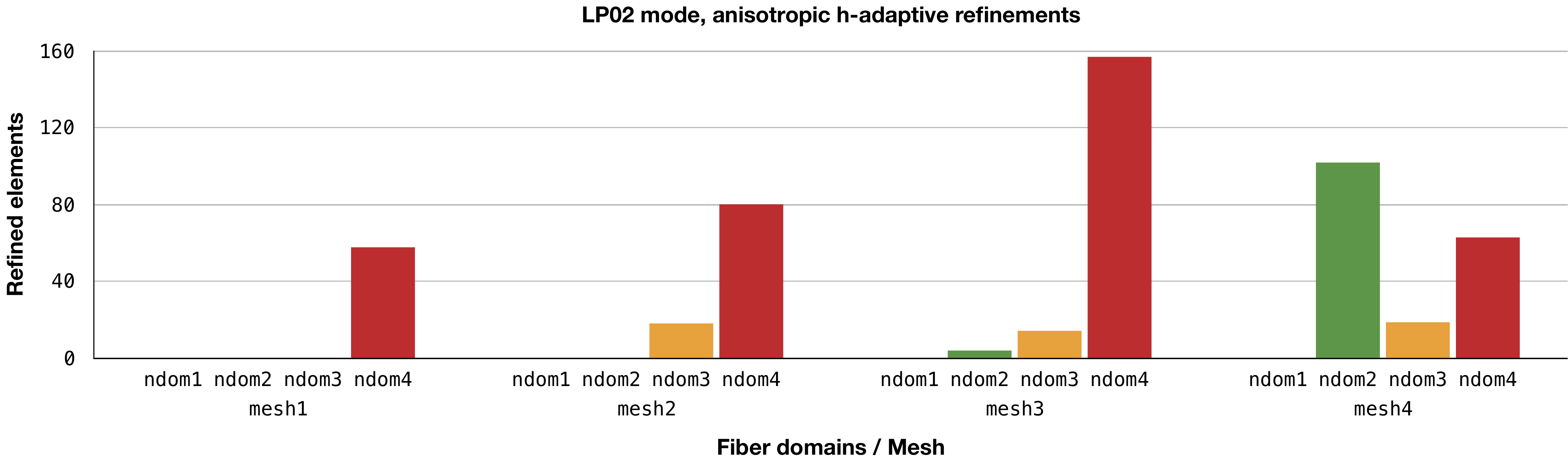}
	\end{subfigure}
	\caption{Radial (anisotropic) $h$-adaptive refinements. For higher-order modes, anisotropic refinements are more computationally efficient to capture the transverse mode profile.}
	\label{fig:aniso-adaptive}
\end{figure}

Next, we repeat the experiment with anisotropic (radial) $h$-refinements. Radial refinements are of interest because higher-order modes are only more oscillatory in the transverse field, but they are not more oscillatory in the direction of propagation. In other words, the guided modes have very similar propagation constants (in fact, higher-order modes oscillate slightly slower than the fundamental mode). Therefore, if the numerical pollution is low for the fundamental mode, we may assume that the resolution in the direction of propagation is ``good enough'' for approximating any higher-order guided modes. Then, radial refinements (in $h$ or $p$) are the more economical way of capturing these modes. For anisotropic $h$-adaptive refinements, depicted in Fig.\ \ref{fig:aniso-adaptive}, a similar pattern emerges for the higher-order modes but the picture is quite different for the $\LP_{01}$ mode. The fundamental mode repeatedly refines elements in the same domain, indicating that the anisotropic refinements do not decrease the local residuals in a way the isotropic ones did. This indicates that the fundamental mode is already well approximated in the transverse field and the residual demands refinements in $z$-direction for better accuracy of the numerical solution. This interplay between the resolution in different directions is critical when studying the pollution error for guided modes.

\begin{figure}[htb]
	\centering
	\begin{subfigure}[b]{0.4\textwidth}
		\includegraphics[width=\textwidth]{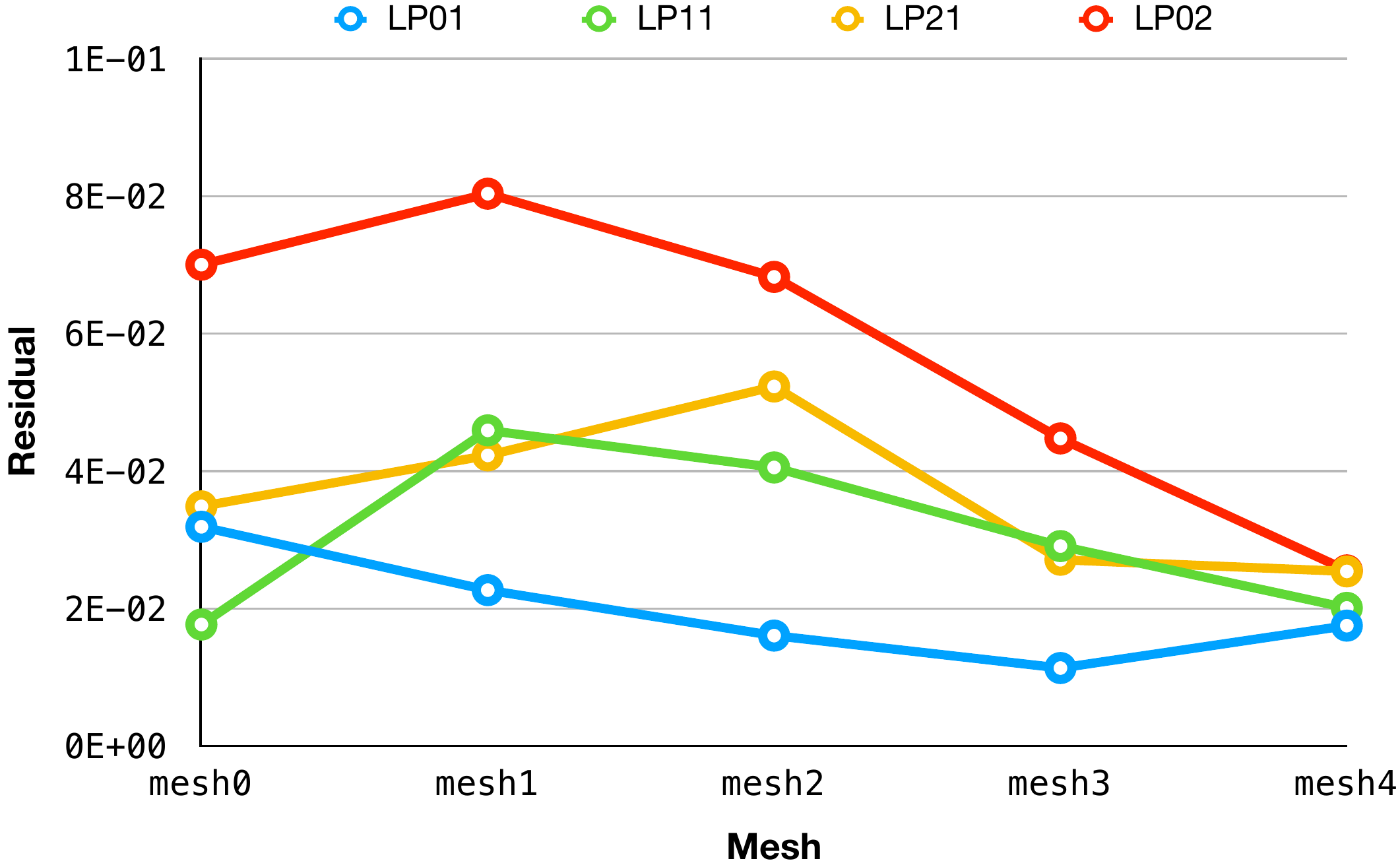}
		\caption{Isotropic $h$-adaptive refinements}
		\label{fig:residual-adapt-iso}
	\end{subfigure}
	\begin{subfigure}[b]{0.4\textwidth}
		\includegraphics[width=\textwidth]{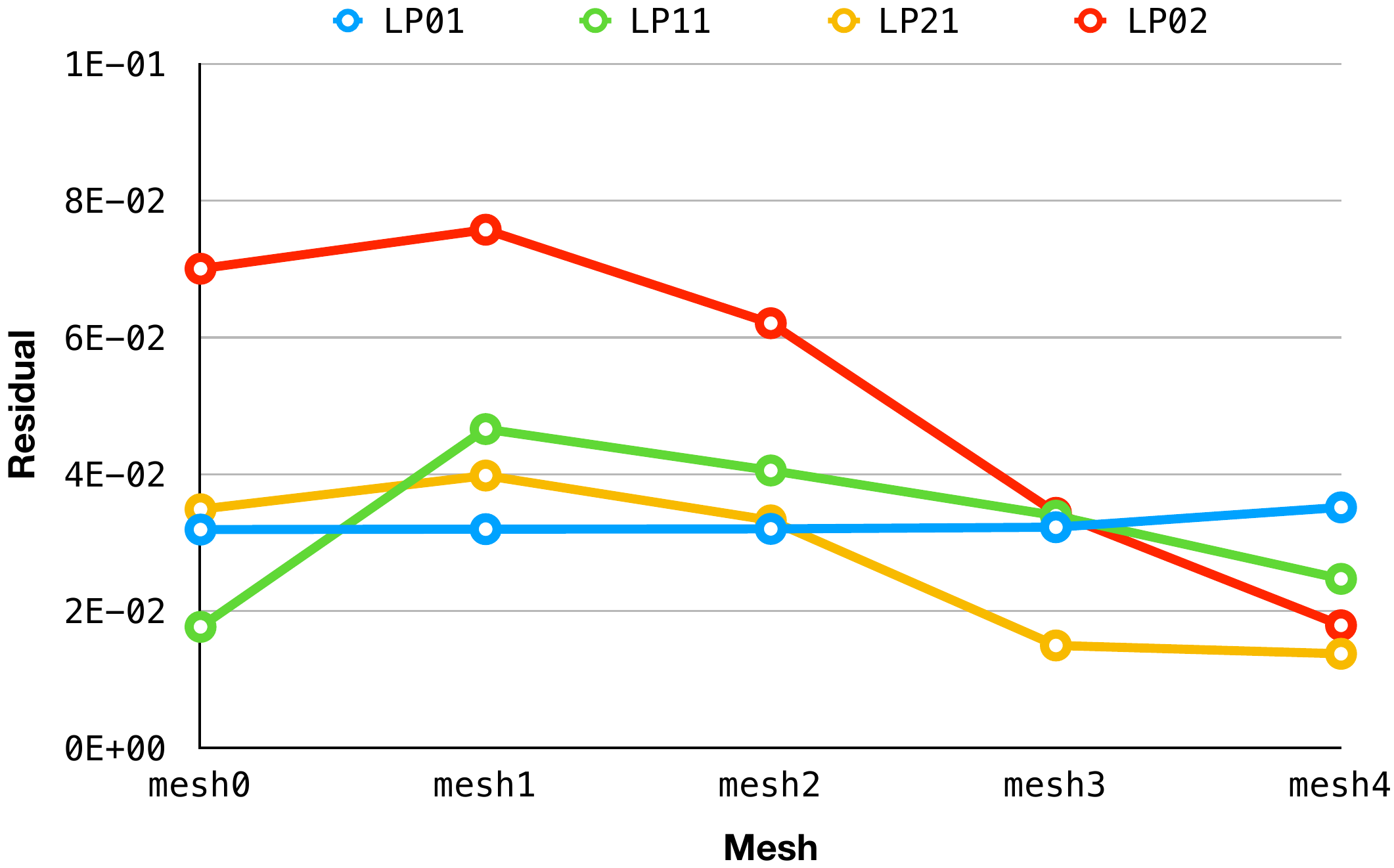}
		\caption{Anisotropic $h$-adaptive refinements}
		\label{fig:residual-adapt-aniso}
	\end{subfigure}
	\caption{Evolution of the DPG residual in adaptive mesh refinements. With our choice of initial geometry discretization, the anisotropic refinements decrease the residual for higher-order modes but not for the fundamental mode. This illustrates that the interplay of the resolution between the transverse direction and the direction of propagation is critical, and the optimal refinement strategy depends on the propagating modes.}
	\label{fig:residual-adapt}
\end{figure}

Fig.~\ref{fig:residual-adapt} shows how the total residual evolves in both scenarios: we observe that the higher-order modes benefit much from anisotropic refinements, making this the preferred choice for improving the numerical solution with fewer degrees of freedom. For the fundamental mode, we find that the residual does not further decrease through anisotropic refinements indicating the mode is captured quite well by the initial current geometry.

\FloatBarrier

\subsection{Load balancing}
\begin{figure}[htb]
	\centering
	\includegraphics[width=0.8\textwidth]{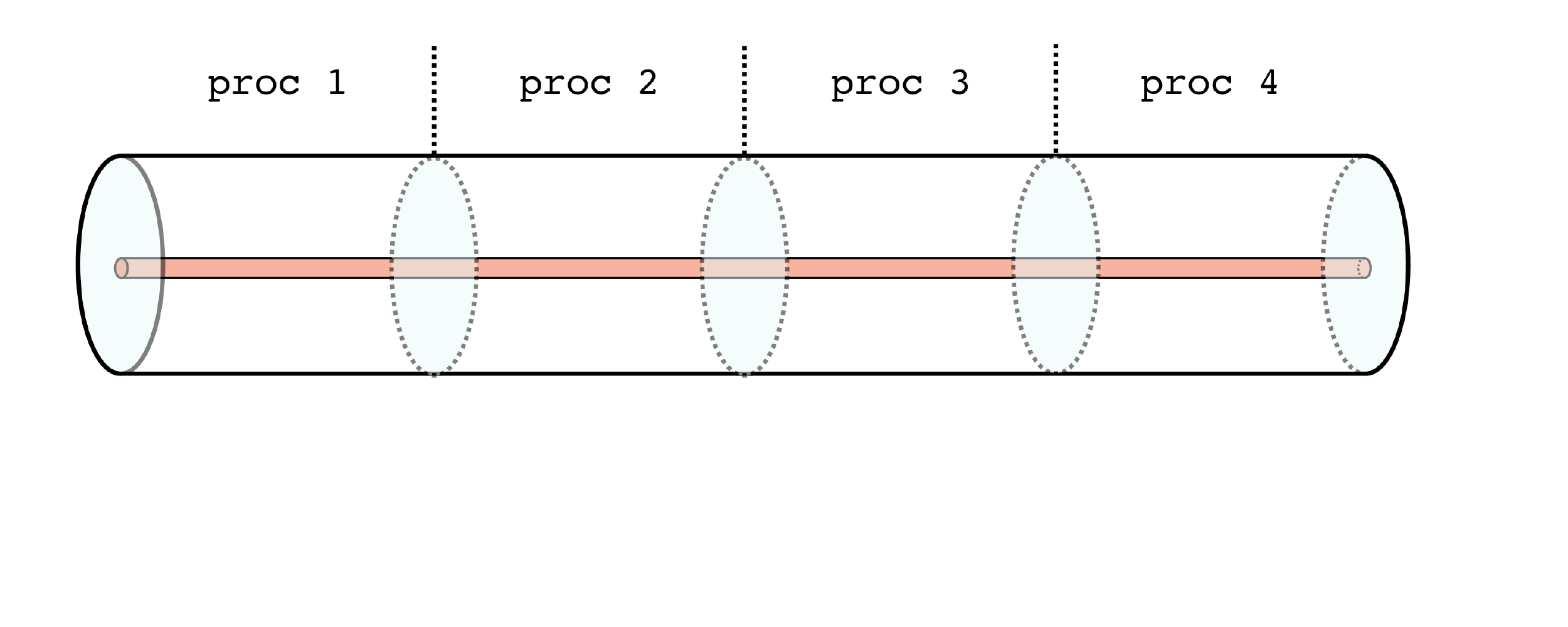}
	\vskip -50pt
	\caption{Initial (static) load distribution in the step-index fiber. An efficient and balanced distribution is achieved by defining subdomains through orthogonal cuts to the fiber axis.}
	\label{fig:LB7-fiber-3}
\end{figure}

In the parallel computation of the fiber problem, we partition the geometry into subdomains, each owned by one distinct MPI process. The rank of each MPI process is the ID of the subdomain it owns. Initially, we partition the fiber directly based on geometric cuts orthogonal to the fiber axis. Fig.~\ref{fig:LB7-fiber-3} illustrates what the partitioning looks like for four subdomains.

This initial static partitioning is a good choice because it keeps the interfaces between subdomains small, making it possible to compute large fibers in parallel with a nested dissection solve and obtain good weak scaling. As the adaptive mesh refinements proceed, we must dynamically repartition the domain to retain load balance. A large number of different repartitioners are available in open source software packages, such as Zoltan \cite{ZoltanOverviewArticle2002}. We believe that in this instance, graph partitioning that strives for minimum cuts is a good choice because it keeps the subdomain interfaces relatively small. In the broken ultraweak Maxwell formulation, we are using element interior dofs (electromagnetic fields: $L^2$) as weights for graph vertices, and trace DOFs on faces (electromagnetic fluxes: $H(\text{curl})$ trace) as weights for the graph's edges. We are omitting connectivities from edge degrees of freedom to provide a sparser graph and accelerate partitioning. ParMetis or PT-Scotch can be used for approximating the partitioning problem. As an alternative, we use a custom dynamic fiber repartitioner that forces orthogonal cuts through the domain while trying to maximize load balance and minimize data migration, similar to recursive coordinate bisection partitioners. This custom repartitioner can perform orders of magnitude faster than graph partitioning because it relies primarily on geometry information.

\begin{figure}[htb]
	\centering
		\includegraphics[width=0.8\textwidth]{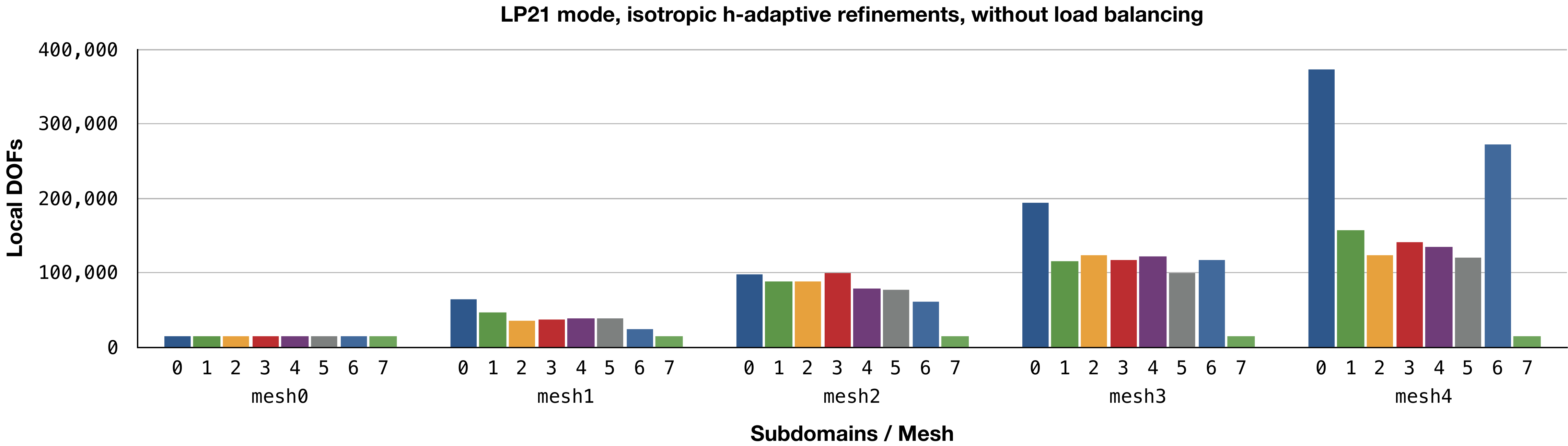}
	\caption{Isotropic $h$-adaptive refinements. The load imbalance increases with every adaptive refinement step. The most refinements are happening at the fiber input (subdomain 0) and within the first few wavelengths of the PML region (subdomain 6).}
	\label{fig:iso-no-LB}
\end{figure}

\begin{figure}[htb]
	\centering
		\includegraphics[width=0.8\textwidth]{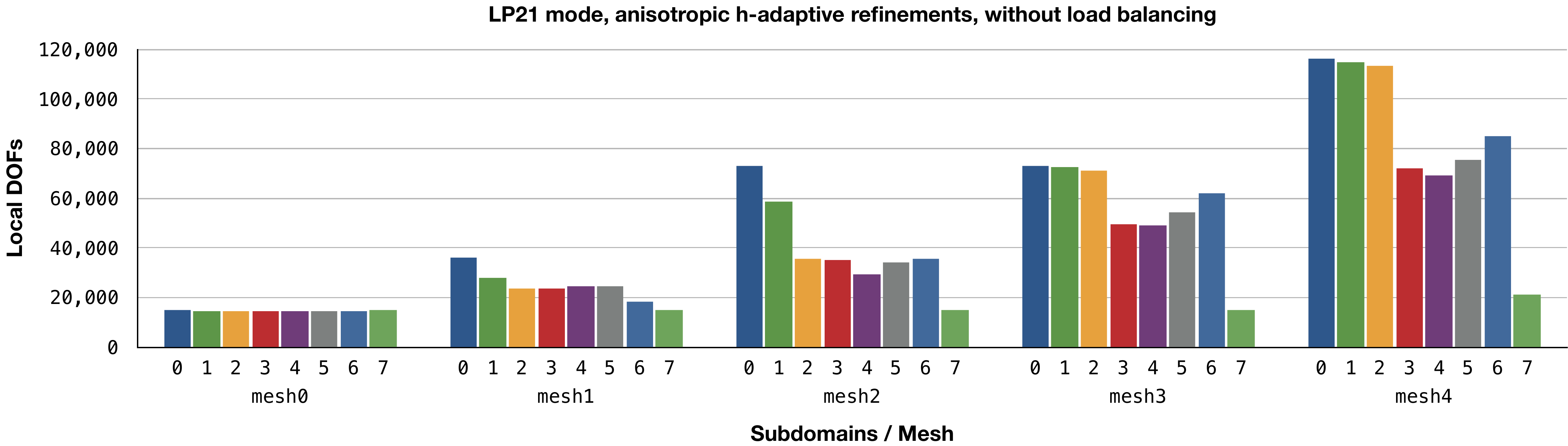}
	\caption{Radial (anisotropic) $h$-adaptive refinements. Load imbalance also occurs with the anisotropic adaptivity but it is less pronounced than in the isotropic case.}
	\label{fig:aniso-no-LB}
\end{figure}

We are studying how the workload in different subdomains changes without repartitioning. Fig.\ \ref{fig:iso-no-LB} and Fig.\ \ref{fig:aniso-no-LB} show the workload per MPI rank in a fiber of 16 wavelengths, partitioned into 8 subdomains, with $h$-adaptive isotropic and anisotropic refinements, respectively. Both plots show the results for the higher-order mode $\LP_{21}$. Here, the workload is simply shown as the number of subdomain interior DOFs, i.e., all solution DOFs that are part of a subdomain excluding the trace DOFs on the subdomain interfaces.


Isotropic refinements lead expectedly to higher load imbalance because each isotropic $h$-refinement increases the local DOFs by another factor of two compared to the anisotropic $h$-refinement. Two observations stand out: firstly, the subdomains closer to the fiber input appear to exhibit higher residuals hence more refinements are observed in that region; secondly, towards the end of the fiber, many refinements are picked up in the sixth subdomain, and almost none in the very last one. The latter observation is an effect from the PML boundary layer at the fiber end. In this short fiber, the PML is active in the last two subdomains (i.e., the layer encompasses about four wavelengths). When the wave enters the boundary layer, it exhibits exponential decay due to the coordinate stretching. This initial decay must be captured accurately by the numerical solution. We see that the DPG residual recognizes the need for more refinements in this region and marks elements in the sixth subdomain. By the time the wave enters the last subdomain, owned by rank seven, it has decayed so far that the residual remains fairly small and almost none of the elements are marked in the adaptive procedure.

\begin{figure}[htb]
	\centering
	\begin{subfigure}[b]{0.8\textwidth}
		\includegraphics[width=\textwidth]{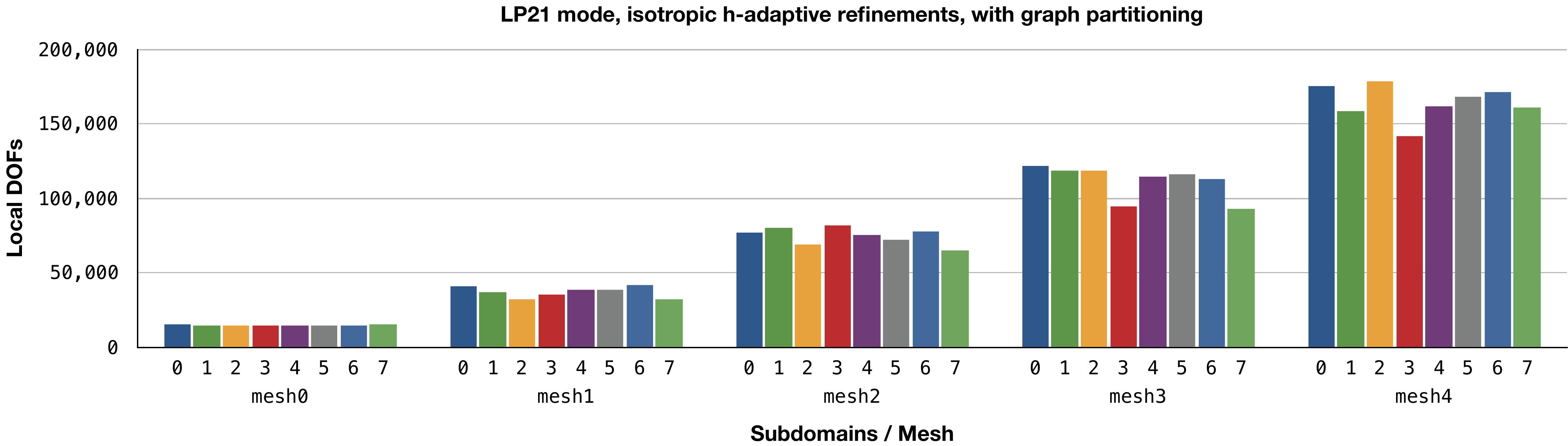}
	\end{subfigure}
	\begin{subfigure}[b]{0.8\textwidth}
		\includegraphics[width=\textwidth]{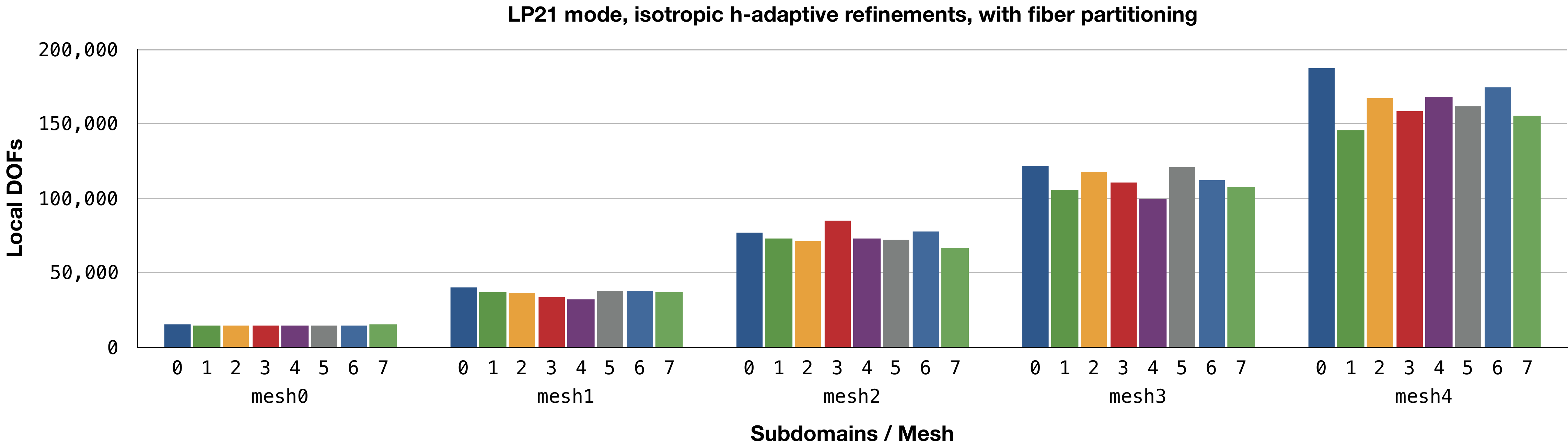}
	\end{subfigure}
	\caption{Load balancing for isotropic $h$-adaptive refinements. Both the graph partitioning and the custom fiber partitioning distribute the workload evenly among the processors in every refinement step.}
	\label{fig:iso-LB}
\end{figure}

It is evident that dynamic load balancing is necessary for computational efficiency in the simulation of the TMI phenomenon or other applications with energy transfer between guided modes. By repartitioning the fiber domain, we obtain significant speedup in the total computation time. Fig.\ \ref{fig:iso-LB} shows the workload per MPI process when the mesh is repartitioned after every isotropic $h$-adaptive refinement step. Both repartitioners, the graph partitioner based on ParMetis as well as the custom fiber partitioner, distribute the workload evenly among the processors. The distribution is not exactly even in either case because both also aim to minimize the size of the interface problem that separates the subdomains. That is, the graph partitioner approximates minimum cuts while the custom fiber partitioner keeps all cuts orthogonal to the fiber axis.

\begin{figure}[htb]
	\centering
	\begin{subfigure}[b]{0.42\textwidth}
		\includegraphics[width=\textwidth]{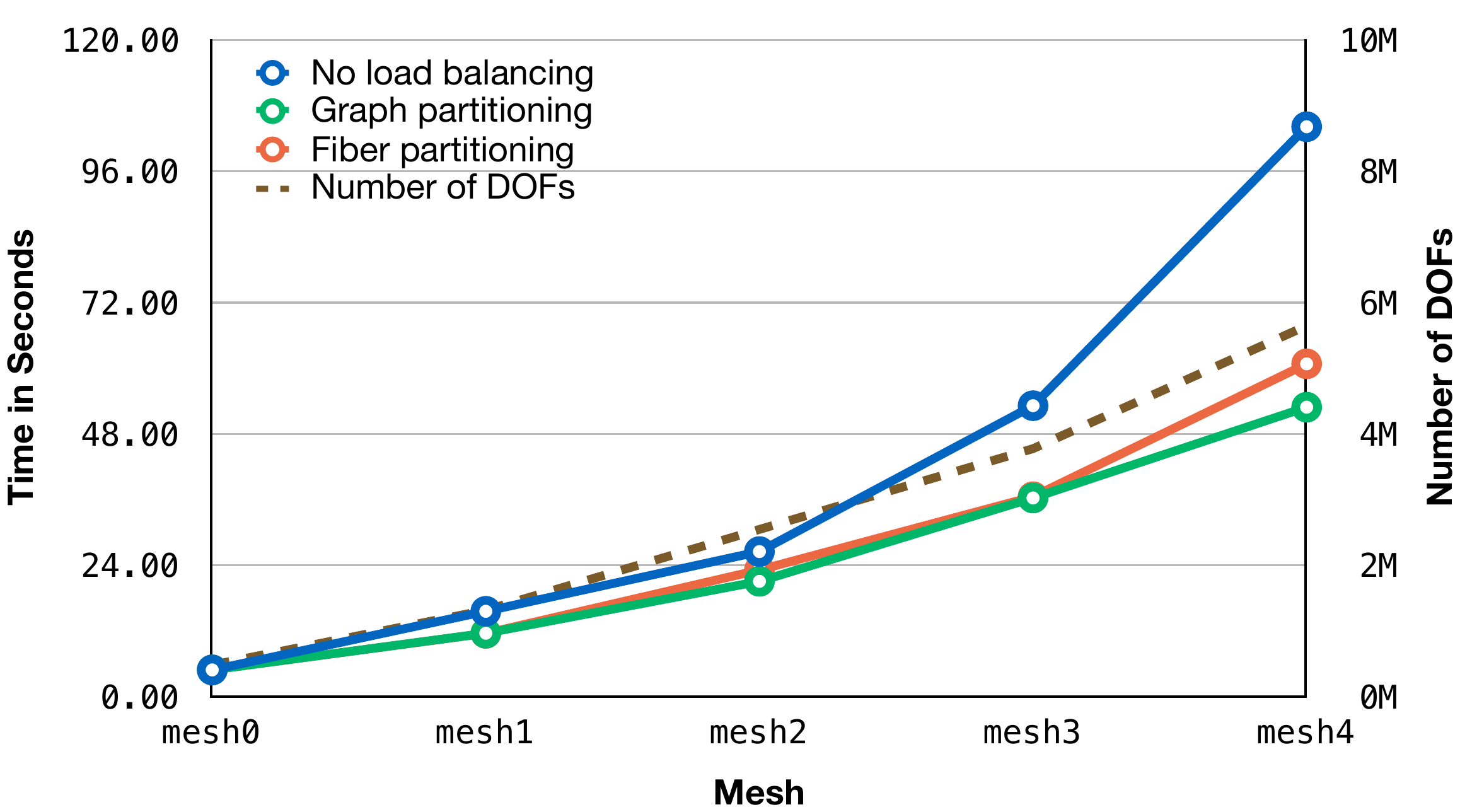}
		\caption{Assembly time}
		\label{fig:assembly-time}
	\end{subfigure}
	\begin{subfigure}[b]{0.39\textwidth}
		\includegraphics[width=\textwidth]{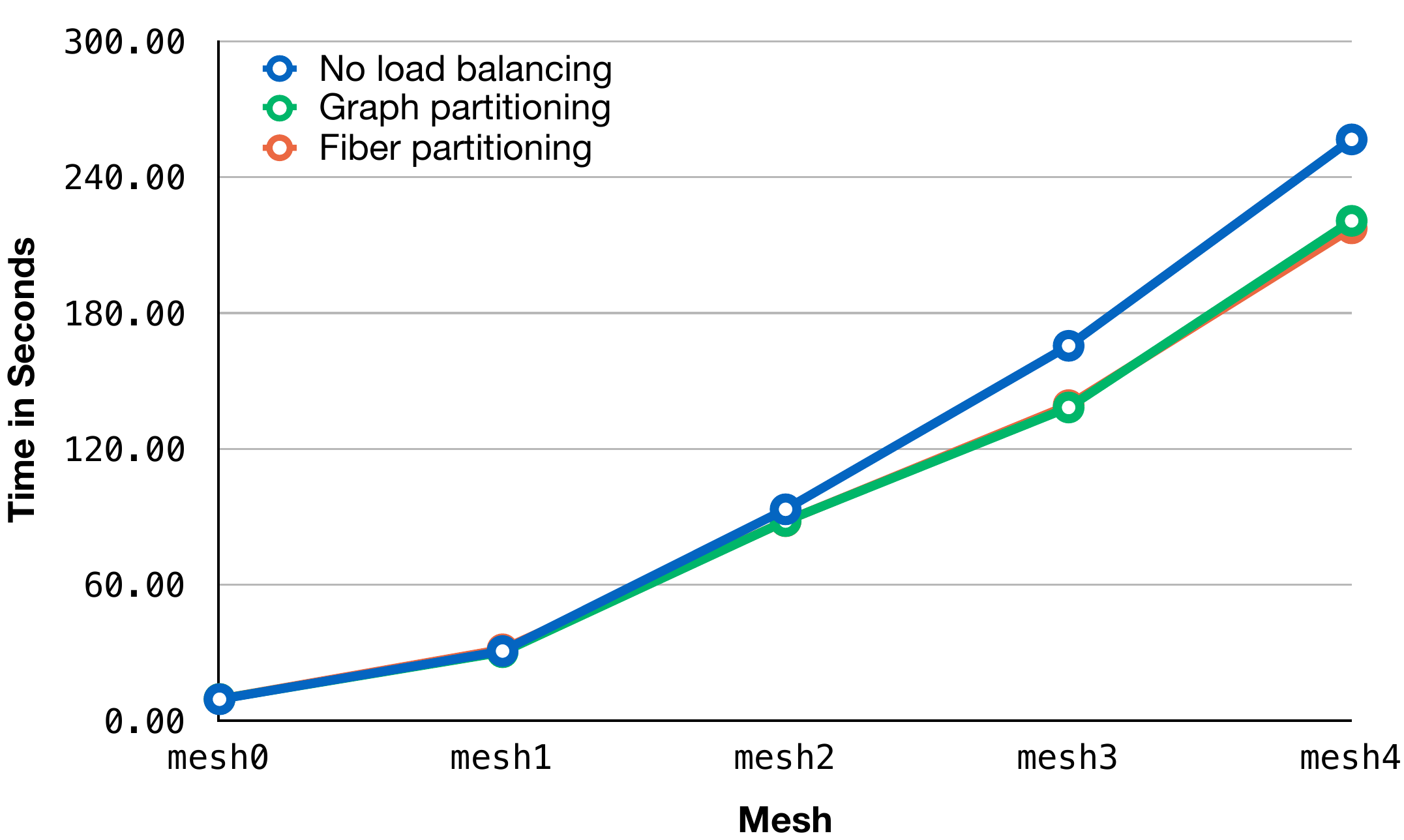}
		\caption{Solve time}
		\label{fig:solve-time}
	\end{subfigure}
	\caption{Computation time with isotropic $h$-adaptive refinements for the guided $\LP_{21}$ mode. The finite element assembly time (including time for element integration) is proportional to the (maximal) workload per subdomain. With load balancing, the assembly time scales corresponding to the increase in the total number of degrees of freedom. The solve time includes the time for analysis, factorization, and linear solve performed by the distributed MUMPS solver. Here, we observe a modest reduction in the computation time when load balancing is performed.}
	\label{fig:assembly-solve-time}
\end{figure}

Fig.~\ref{fig:assembly-time} displays the computation time for the distributed finite element assembly for each of the isotropic $h$-adaptive meshes. The assembly time includes the time for element integration and the assembly of every sparse subdomain stiffness matrix and load vector. In DPG methods, the assembly can be a substantial part of the entire time to solution because of the computation of optimal test functions in the enriched test space. On the other hand, it is conveniently parallel and will exhibit good parallel scaling as long as the workload is balanced. We observe that with load balancing the increase in computation time corresponds to the total increase in the number of DOFs, as expected; in the imbalanced case, where some processes finish early and remain idle until the MPI process with the maximal workload is done, the assembly time increases unproportionally.

The parallel solve time is depicted in Fig.~\ref{fig:solve-time}; it shows the total time for analysis, factorization, and linear solve (forward and backward elimination) performed by the distributed MUMPS solver \cite{amestoy2001mumps}. Both load balancing strategies result in a modest reduction of the solve time and perform equally well. We expect that larger gains could be obtained for parallel nested dissection solves on large problem instances, but we have not yet explored this sufficiently. Of course, the time spent on load balancing, i.e., partitioning and data migration, should be taken into account, as well. However, for the problem size considered here, we have found that the computation time for these tasks is almost negligibly small (less than one second per mesh).

\FloatBarrier

\section{Summary}
%
%

We have studied the effect of numerical pollution for guided wave problems with many wavelengths in conforming DPG finite element discretizations. We found that the pollution error in the ultraweak DPG setting is primarily present in the form of a diffusive effect that causes attenuation of the propagating wave. This is consistent with previous observations for the DPG method where the phase error appears to be relatively small. Furthermore, we were able to corroborate theoretical estimates by Melenk and Sauter, indicating a logarithmic dependency of the pollution error on the polynomial order of approximation. Based on our numerical results, we agree that the best strategy for resolving many wavelengths is one based on $hp$-refinements that first ``capture'' the wave (discretizing each wavelength by a few elements with moderate $p$ so that $\omega h/p$ is sufficiently small) and subsequently increase the polynomial order $p$ with $\mathcal{O}(\log \omega)$ if needed while keeping $\omega h$ constant. For the 3D waveguide problem where the wave is truly a superposition of transverse modes propagating in one direction, we emphasize that, perhaps counterintuitively, controlling the pollution error asymptotically required additional degrees of freedom in both the transverse direction and the direction of propagation.

For controlling the error in a multi-mode waveguide, we described suitable adaptive refinement strategies. The error indicator we used is based on the Riesz representation of the residual in the ultraweak DPG formulation. With broken test spaces, this indicator is computed locally and in parallel. We have shown the importance of adaptive refinements in capturing different propagating modes in a fiber waveguide. For such problems, the efficacy of the DPG residual is remarkable as it sharply recognizes where the ``energy'' of the solution is located. The numerical results in the weakly-guiding optical fiber waveguide indicate that anisotropic $h$-adaptive refinements can be a much more efficient choice than isotropic $h$-adaptive refinements for these problems.

In a distributed-memory computation, the repartitioning of the multi-mode fiber becomes essential to maintain load balance. We observed that the PML boundary layer requires additional refinements especially as the wave enters the PML region. Through dynamic repartitioning, the load imbalance can be nearly eliminated and the time to solution decreases significantly. This is especially important in problems with localized features in the solution such as the transverse mode instability in fiber amplifiers.

\section*{Acknowledgement}
This work was partially supported by AFOSR FA9550-17-1-0090 and AFOSR 18RDCOR018.

\bibliographystyle{abbrv}

\bibliography{./shortref,./ref,./ref_laser}

\end{document}